# Nonlinear emergent macroscale PDEs, with error bound, for nonlinear microscale systems

J. E. Bunder [*]     A. J. Roberts [†]

2018-06-27

## Abstract

Many physical systems are formulated on domains which are relatively large in some directions but relatively thin in other directions. We expect such systems to have emergent structures that vary slowly over the large dimensions. Common mathematical approximations for determining the emergent dynamics often rely on self-consistency arguments or limits as the aspect ratio of the 'large' and 'thin' dimensions becomes nonphysically infinite. Here we extend to nonlinear dynamics a new approach [IMA J. Appl. Maths, doi:10.1093/imamat/hxx021] which analyses the dynamics at each cross-section of the domain via a rigorous multivariate Taylor series. Then centre manifold theory supports the global modelling of the system's emergent dynamics in the large but finite domain. Interactions between the cross-section coupling and both fast and slow dynamics determines quantitative error bounds for the nonlinear modelling. We illustrate the methodology by deriving the large-scale dynamics of a thin liquid film, where the film is subject to a Coriolis force induced by a rotating substrate. The approach developed here quantifies the accuracy of known approximations, extends such approximations to mixed order modelling, and may open previously intractable modelling issues to new tools and insights.

## Contents

[*]School of Mathematical Sciences, University of Adelaide, South Australia. mailto:judith.bunder@adelaide.edu.au , http://orcid.org/0000-0001-5355-2288

[†]School of Mathematical Sciences, University of Adelaide, South Australia. mailto:anthony.roberts@adelaide.edu.au , http://orcid.org/0000-0001-8930-1552

# 1   Introduction

Many systems of interest in science and engineering occur in a domain with disparate length scales (Davis 2017, e.g.): often a fine structure is modulated on a much larger scale (Mielke 1992, e.g.). Such disparate scales often are a major challenge in computational simulations (Rüde et al. 2016, p.14, e.g.). Two classic examples are Taylor–Couette flow (Iooss and Adelmeyer 1992, e.g.) and Benard convection (Segel 1969, e.g.). Often the fine-scale detail is crucial to the accurate modelling of a multiscale system, but with multiple length scales a simulation resolving the physical fine scale is not only prohibitively inefficient and severely constrained by memory limitations, it is also an arduous task to analyse simulation data generated at a scale much smaller than the scale of interest. This article develops a unified mathematical theory to use multiple length scales to reduce the full set of nonlinear governing equations to a simplified evolution equation, with quantified error, enabling more efficient simulations and analysis. This theoretical methodology should be able to better justify and illuminate many extant longwave and homogenisation theories (Bakhvalov and Panasenko 1989; Cross and Hohenberg 1993; Oron, Davis, and Bankoff 1997; Roberts and Li 2006; Pavliotis and Stuart 2008).

    We focus on a class of multiscale systems whose physical domain is 'large' in multiple directions, but have a relatively 'thin' cross-section in the other directions. As an specific example, Section 1.1 takes as prescribed a variant of the integrated boundary layer PDEs for a thin liquid film of Newtonian fluid spreading over a planar rotating surface (Chang 1994), and rigorously derives a simpler lubrication model of the nonlinear flow of the film (Wilson, Hunt, and Duffy 2000; Oron, Davis, and Bankoff 1997, §II.K, e.g.). Appendix A lists computer algebra code for deriving this lubrication model, with the code written to be readily adaptable to a wide range of systems. Thin liquid flows are important in biology, physics, and engineering, as well as in the environment. They may be of common liquids such as water or oil, or of more rheologically complex fluids, and display interesting nonlinear wave patterns (Oron, Davis, and Bankoff 1997; Wilson, Hunt, and Duffy 2000; Murisic et al. 2011; Lam et al. 2015). Other examples of systems amenable to our methodology include flood and tsunami modelling (Noakes, King, and Riley 2006; Bedient and Huber 1988; LeVeque, George, and Berger 2011), pattern formation (Newell and Whitehead 1969; Cross and Hohenberg 1993; Westra, Binks, and Water 2003), wave interactions (Nayfeh and Hassan 1971; Griffiths, Grimshaw, and Khusnutdinova 2006), elastic shells (Naghdi 1973; Mielke 1986; Lall, Krysl, and Marsden 2003), and microstructured materials (Romanazzi, Bruna, and Howey 2016).

    Section 2 defines the generic nonlinear PDE system to which the methodology applies, and defines the 'large' but 'thin' multiscale domain on which the system evolves. In pattern evolution problems a 'thin' domain variable is the phase of the underlying small scale pattern (Roberts 2015, §3.3, e.g.). A first step in the reduction of such PDEs was taken by Mielke (1992), but the analysis required solutions to exist for all $\mathbb{R}$-time and for all $\mathbb{R}^m$-space which precludes initial/boundary value problems. Here we analyse the dynamics of a general PDE in a thin

cross-section of the large domain by constructing a multivariate Taylor expansion for the local spatial structures and analysing the evolution of the coefficients. Being the union of local-space-time modelling means the approach is valid everywhere outside of boundary layers (Roberts 1992, e.g.) and initial transients (Roberts 1989, e.g.). Section 3 details how to capture the emergent behaviour of the system at every chosen cross-section via constructing a set of generalised eigenvectors which span the centre subspace on which the slow system evolves. Based upon these eigenvectors the system's centre manifold is parametrised and an emergent nonlinear PDE derived for the slow evolution. The order of the multivariate Taylor expansion determines the order of accuracy of the derived slow evolution, and Lagrange's Remainder Theorem provides a novel exact expression for the error of this slow evolution.

This article extends the methodology initially developed by Roberts (2015) to derive the long, slow space evolution of nonlinear PDE systems on a long one-dimensional physical domain with a relatively thin cross-section. Roberts and Bunder (2017) then developed the methodology to linear systems that have a domain that is large in multiple dimensions. Here we further develop the methodology to nonlinear PDEs with multiple large dimensions. Our methodology is different to many other methods which derive the large-scale evolution in that here no asymptotic limit is required for the scale separation between the large domain and the thin cross-section (Kondic 2003; Náraigh and Thiffeault 2010, e.g.). Specifically, we have no requirement that a scale separation parameter (often named $\epsilon$) must be asymptotically small; we only require that such a small-large scale separation exists so that we can establish a centre-stable subspace (Assumption 3), and then our approach is valid at finite scale separation.

## 1.1 Example of a rotating shallow fluid flow

As an example application of some of the results, consider the flow of a shallow layer of fluid on a solid flat rotating substrate, such as in spin coating (Wilson, Hunt, and Duffy 2000; Oron, Davis, and Bankoff 1997, §II.K, e.g.) or large-scale shallow water waves (Dellar and Salmon 2005; Hereman 2009, e.g.). Let $\vec{x} = (x_1, x_2)$ parametrise location on the rotating substrate, and let the fluid layer have thickness $h(\vec{x}, t)$ and move with depth-averaged horizontal velocity $\vec{v}(\vec{x}, t) = (v_1, v_2)$. We take as given (with its simplified physics) that the governing set of PDEs is the nonlinear system

$$\frac{\partial h}{\partial t} = -\nabla \cdot (h\vec{v}), \tag{1a}$$

$$\frac{\partial \vec{v}}{\partial t} = \begin{bmatrix} -b & f \\ -f & -b \end{bmatrix} \vec{v} - (\vec{v} \cdot \nabla)\vec{v} - g\nabla h, \tag{1b}$$

where $b$ represents viscous bed drag, $f$ is the Coriolis coefficient, $g$ is the acceleration due to gravity, and we neglect surface tension. The PDEs (1) are similar to those used by Dellar and Salmon (2005, eq. (21)), but with only one component of the Coriolis force and the addition of viscous drag, and also similar to that used by Hereman (2009, eqs. (22)–(24)), but with a flat substrate.

For such a shallow fluid flow, the horizontal gradient $\nabla$ of quantities are mostly relatively small (Davis 2017, e.g.). Then the flow driven by variations of film thickness, $\nabla h$, is approximately balanced in (1b) by the rotation and the substrate drag, leading to the velocity field

$$\vec{v} \approx \frac{-g}{b^2 + f^2} \begin{bmatrix} b & f \\ -f & b \end{bmatrix} \nabla h.$$

Substituting this balance in the conservation of mass equation (1a) derives the single component, 'lubrication', model

$$\frac{\partial h}{\partial t} \approx \frac{gb}{b^2 + f^2} \nabla \cdot (h \nabla h). \qquad (2)$$

Having just one component, we can use the PDE model (2) as a simpler description of the shallow fluid dynamics. But PDE (2) is an approximation to the 'original' PDE (1), and so three outstanding questions are: can we find a rigorous error? can the analysis be extended to higher order? and can such an approach apply generally? Our answer to all three questions is yes.

Returning to the original system of PDEs (1) and defining the system field $u(\vec{x}, t) = (h, v_1, v_2)$, we rewrite system (1) as one equation, while also combining terms of like order:

$$\frac{\partial u}{\partial t} = \begin{bmatrix} 0 & 0 & 0 \\ 0 & -b & f \\ 0 & -f & -b \end{bmatrix} u + \begin{bmatrix} 0 & -1 & 0 \\ -g & 0 & 0 \\ 0 & 0 & 0 \end{bmatrix} \partial_{x_1} u + \begin{bmatrix} 0 & 0 & -1 \\ 0 & 0 & 0 \\ -g & 0 & 0 \end{bmatrix} \partial_{x_2} u \qquad (3a)$$

$$+ \begin{bmatrix} (\partial_{x_1} u)^\mathsf{T} \mathfrak{M}_{10,0} u \\ (\partial_{x_1} u)^\mathsf{T} \mathfrak{M}_{10,1} u \\ (\partial_{x_1} u)^\mathsf{T} \mathfrak{M}_{10,2} u \end{bmatrix} + \begin{bmatrix} (\partial_{x_2} u)^\mathsf{T} \mathfrak{M}_{01,0} u \\ (\partial_{x_2} u)^\mathsf{T} \mathfrak{M}_{01,1} u \\ (\partial_{x_2} u)^\mathsf{T} \mathfrak{M}_{01,2} u \end{bmatrix}, \qquad (3b)$$

where matrices

$$\mathfrak{M}_{10,0} = \begin{bmatrix} 0 & -1 & 0 \\ -1 & 0 & 0 \\ 0 & 0 & 0 \end{bmatrix}, \quad \mathfrak{M}_{10,1} = \begin{bmatrix} 0 & 0 & 0 \\ 0 & -1 & 0 \\ 0 & 0 & 0 \end{bmatrix}, \quad \mathfrak{M}_{10,2} = \begin{bmatrix} 0 & 0 & 0 \\ 0 & 0 & 0 \\ 0 & -1 & 0 \end{bmatrix},$$

$$\mathfrak{M}_{01,0} = \begin{bmatrix} 0 & 0 & -1 \\ 0 & 0 & 0 \\ -1 & 0 & 0 \end{bmatrix}, \quad \mathfrak{M}_{01,1} = \begin{bmatrix} 0 & 0 & 0 \\ 0 & 0 & -1 \\ 0 & 0 & 0 \end{bmatrix}, \quad \mathfrak{M}_{01,2} = \begin{bmatrix} 0 & 0 & 0 \\ 0 & 0 & 0 \\ 0 & 0 & -1 \end{bmatrix}. \qquad (3c)$$

The system's first line (3a) is the linear part, whereas the second line (3b) contains the nonlinear terms which encode inertial acceleration. When the nonlinear terms (3b) are negligible, such as for very viscous fluids with low Reynolds numbers, the method described by Roberts and Bunder (2017) derives a slow *linear* PDE approximation of (3), but the analysis and results of that article do not account for nonlinear dynamics. Herein we develop the methodology and theory to account for nonlinear effects.

For the shallow fluid flow described by (3), the 'large' domain is some 'physical' subset of the $x_1 x_2$-plane, and the 'thin' cross-section is the three-components of $u = (h, v_1, v_2)$. The aim is to capture the long, slow behaviour of the original $u(\vec{x}, t)$ field in a one-component slow field $U(\vec{x}, t)$ (instead of the three components which describe the thin cross-section), and to construct a PDE for $U(\vec{x}, t)$ which is correct to some order $N$ in spatial derivatives, with known error. In general, higher orders $N$ are potentially able to capture more extreme spatial fluctuations but may not be structurally stable, and so we address up to some low–moderate order $N$. This restricts $U(\vec{x}, t)$ to describing long, relatively gradual, spatial variations of the original field $u(\vec{x}, t)$. Section 3.2 defines the slow field $U(\vec{x}, t)$ in the general case and proves that it describes the behaviour of the original microscale field on the slowly evolving centre manifold (or slow manifold, as is the case in this shallow fluid example).

To determine the nature of the slow field $U$ we first seek to understand the 'slow' and 'fast' evolution of the lowest-order linearisation of the system field $u$ (Assumption 3 elaborates the general case). The linear dynamics of $u$ are dominantly characterised by the lowest order linear term in (3), the term $\mathfrak{L}_{\vec{0}} u$ where here
$$\mathfrak{L}_{\vec{0}} = \begin{bmatrix} 0 & 0 & 0 \\ 0 & -b & f \\ 0 & -f & -b \end{bmatrix}.$$

The eigenvalues of $\mathfrak{L}_{\vec{0}}$ indicate that $u$ evolves on a one-dimensional slow subspace (one zero eigenvalue) and a two-dimensional stable subspace (two eigenvalues, $-b \pm if$, with negative real part). The stable part of $u$ decays relatively quickly, roughly like $e^{-bt}$, while the slow component of $u$, namely $h$, evolves on the one-dimensional slow subspace. In the notation introduced by Assumption 3 and Section 3, there exists a slow subspace of dimension $m = 1$ with right and left eigenvectors $V^{\vec{0}} = (1, 0, 0)$ and $Z^{\vec{0}} = (1, 0, 0)$, and eigenvalue $A_{\vec{0}} = 0$.

Once the lowest-order linear dynamics of the system field $u$ are known from matrix $\mathfrak{L}_{\vec{0}}$, we construct the generalised eigenvectors $\tilde{\mathcal{V}}^{\vec{n}}$, for $|\vec{n}| \leqslant N$, which span the spatially-local slow subspace of the linear system (3a) to the specified order of accuracy $N$. This order of accuracy $N$ is that of a local multivariate Taylor expansion of the field $u(\vec{x}, t)$. The advantages of such a Taylor expansion are that not only does it provide a straightforward way to increase the order of accuracy $N$, it also provides a rigorous error term from Lagrange's Remainder Theorem. Section 3.1 discusses the general construction of the eigenvectors $\tilde{\mathcal{V}}^{\vec{n}}$, which in essence detail local out-of-equilibrium structures, and then Section 3.2 models the full nonlinear system to derive the slow manifold PDE.

Appendix A lists computer algebra Reduce[1] code which applies the general theory of Sections 2 and 3 to determine the slow PDE for any order $N$ of the Taylor expansion, thus constructing PDEs which describe the slow $U = h$ field evolution of the shallow fluid dynamics to various orders of accuracy. See Sections 2 and 3 for details that justify (2) for order $N = 2$, and that for order $N = 3$ the slow PDE is the 2D advection-diffusion PDE

$$\frac{\partial h}{\partial t} \approx \frac{bg}{b^2 + g^2} \nabla^2 h + \sum_{|\vec{\ell}_1|=0}^{3} \sum_{|\vec{\ell}_2|=0}^{|\vec{\ell}_1|} a_{\vec{\ell}_1 \vec{\ell}_2} (\partial_{\vec{x}}^{\vec{\ell}_1} h)(\partial_{\vec{x}}^{\vec{\ell}_2} h), \tag{4a}$$

with multi-indices $\vec{\ell}_{1,2} \in \mathbb{N}_0^2$, symmetry $a_{\vec{\ell}_1 \vec{\ell}_2} = a_{\vec{\ell}_2 \vec{\ell}_1}$ when $|\vec{\ell}_1| = |\vec{\ell}_2|$, and nonzero constant coefficients [2]

$$a_{(01)(01)} = a_{(10)(10)} = a_{(02)(00)} = a_{(20)(00)} = \frac{bg}{b^2 + f^2},$$
$$a_{(03)(01)} = a_{(30)(10)} = a_{(12)(10)} = a_{(21)(01)} = -\frac{bg^2(b^2 - f^2)}{(b^2 + f^2)^3},$$
$$a_{(03)(10)} = a_{(21)(10)} = -a_{(30)(01)} = -a_{(12)(01)} = -2f \frac{b^2 g^2}{(b^2 + f^2)^3}. \tag{4b}$$

The four $c$ coefficients equal to $bg/(b^2 + f^2)$ correspond to the $N = 2$ approximation (2), but the rigorous derivation from Appendix A on the slow subspace empowers a far richer description of the slow dynamics in the $x_1x_2$-plane. Further, Section 3.2 provides an exact expression (45) for the error of slow PDEs such as (4a); Appendix A.3 details how components of this error are constructed for this shallow fluid flow.

---

[1] Reduce [http://reduce-algebra.com/] is a free, fast, general purpose, computer algebra system.

[2] Equation (3b) describes two nonlinear terms which, following the derivation in Section 3, result in nonlinear terms with constant coefficients $a^j_{\vec{\ell}_1 \vec{\ell}_2}$ for $j = 1, 2$. But in this example, since the nonlinear terms are both second order ($P_j = 2$ for $j = 1, 2$), to obtain (4a) we set $a^1_{\vec{\ell}_1 \vec{\ell}_2} + a^2_{\vec{\ell}_1 \vec{\ell}_2} = a_{\vec{\ell}_1 \vec{\ell}_2}$.

Executing the computer algebra code in Appendix A to obtain slow PDEs of higher orders is straightforward. Although the computational time increases rapidly with $\mathsf{N}$, in principle the code is applicable to any order $\mathsf{N}$. For example, the slow PDE for order $\mathsf{N}=4$ is

$$\frac{\partial \mathsf{h}}{\partial \mathsf{t}} \approx \frac{\mathsf{bg}}{\mathsf{b}^2+\mathsf{g}^2}\nabla^2\mathsf{h} - \frac{\mathsf{bg}^2(\mathsf{b}^2-\mathsf{f}^2)}{(\mathsf{b}^2+\mathsf{f}^2)^3}\nabla^4\mathsf{h} + \sum_{|\vec{\ell}_1|=0}^{4}\sum_{|\vec{\ell}_2|=0}^{|\vec{\ell}_1|}\mathsf{a}_{\vec{\ell}_1\vec{\ell}_2}(\partial_{\vec{x}}^{\vec{\ell}_1}\mathsf{h})(\partial_{\vec{x}}^{\vec{\ell}_2}\mathsf{h}), \quad (5a)$$

with the same nonzero constant coefficients $\mathsf{a}_{\vec{\ell}_1\vec{\ell}_2}$ given in (4b), as well as

$$\mathsf{a}_{(40)(00)} = \mathsf{a}_{(00)(40)} = \tfrac{1}{2}\mathsf{a}_{(22)(00)} = -\frac{\mathsf{bg}^2(\mathsf{b}^2-\mathsf{f}^2)}{(\mathsf{b}^2+\mathsf{f}^2)^3}. \quad (5b)$$

For small enough damping, $\mathsf{b}<\mathsf{f}$, the fourth order hyperdiffusion in (5a) makes the model structurally unstable. As is often necessary, for $\mathsf{b}<\mathsf{f}$ one would then regularise the model as in the Benjamin, Bona, and Mahony (1972) regularised long wave equation. Notwithstanding such practical regularisation, our derived error expression (45) applies and is useful for as long as the spatial gradients in the solutions to (5a) remain small enough.

The following sections develop theoretical support for the derivation of nonlinear slow PDEs such as (2), (4) and (5), and derive the novel exact algebraic expression (45) for the error.

## 2 Local expansion of general nonlinear dynamics

Consider some multiscale spatial domain $\mathbb{X}\times\mathbb{Y}$ where $\mathbb{X}$ is some open domain of large macroscale extent, and $\mathbb{Y}$ is a 'relatively small' microscale domain (in some Hilbert space). We analyse the dynamics of some field $\mathsf{u}$ within the multiscale spatial domain $\mathbb{X}\times\mathbb{Y}$ and determine the emergent behaviour of this field on the macroscale; that is, we aim to derive a description, over some time interval $\mathbb{T}$, of the long, slow $\mathsf{u}$ field dynamics on the macroscale domain $\mathbb{X}$ while accounting for the fine details in the microscale domain $\mathbb{Y}$ in an 'averaged', 'homogenised' or 'slaved' sense (Roberts 2015; Roberts and Bunder 2017). As the domain $\mathbb{Y}$ is a small cross-section of the full domain of the system, a description of the large-scale behaviour should not involve fluctuations across $\mathbb{Y}$ as dynamic variables.

We consider the field $\mathsf{u}(\vec{x},\mathsf{y},\mathsf{t})$ in a given Hilbert space $\mathbb{U}$ (finite or infinite dimensional), where $\mathsf{u}:\mathbb{X}\times\mathbb{Y}\times\mathbb{T}\to\mathbb{U}$ is a function of $\mathsf{M}$-dimensional position $\vec{x}\in\mathbb{X}\subseteq\mathbb{R}^\mathsf{M}$, cross-sectional position $\mathsf{y}\in\mathbb{Y}$, and time $\mathsf{t}\in\mathbb{T}\subseteq\mathbb{R}$. We suppose the field $\mathsf{u}(\vec{x},\mathsf{y},\mathsf{t})$ satisfies some specified nonlinear PDE in the form

$$\frac{\partial \mathsf{u}}{\partial \mathsf{t}} = \mathfrak{L}[\mathsf{u}] + \mathsf{f}[\mathsf{u}] = \sum_{|\vec{\mathsf{k}}|=0}^{\infty}\mathfrak{L}_{\vec{\mathsf{k}}}\partial_{\vec{x}}^{\vec{\mathsf{k}}}\mathsf{u} + \mathsf{f}[\mathsf{u}], \quad (6)$$

where $\mathsf{f}[\mathsf{u}]:\mathbb{U}\to\mathbb{U}$ is a 'strictly' nonlinear function of field $\mathsf{u}$ and its derivatives, the $\mathfrak{L}_{\vec{\mathsf{k}}}$ are linear operators (in $\mathsf{y}$), the mulrivariate (mixed) derivative

$$\partial_{\vec{x}}^{\vec{\mathsf{k}}} := \frac{\partial^{|\vec{\mathsf{k}}|}}{\partial \mathsf{x}_1^{\mathsf{k}_1}\partial \mathsf{x}_2^{\mathsf{k}_2}\cdots\partial \mathsf{x}_\mathsf{M}^{\mathsf{k}_\mathsf{M}}}$$

is of order $|\vec{k}| = k_1 + k_2 + \cdots + k_M$, and where the infinite sum in the PDE (6) is notionally written as being over all possible multi-indices $\vec{k} \in \mathbb{N}_0^M$ (as usual, the set of natural numbers $\mathbb{N}_0 := \{0, 1, 2, \ldots\}$), although in practice there will be only a finite number of terms in this sum.

In application to fluid or heat convection the nonlinear term is the quadratic $f[u] = \vec{u} \cdot \vec{\nabla} \vec{u}$, whereas for the Ginzburg–Landau equation $f[u] = u^3$. Consequently we consider nonlinearities that are sums of products of such factors.

**Assumption 1.** *The nonlinear function $f[u]$ may be written as, or usefully approximated as, a sum of products of $u$ and its derivatives:*

$$f[u] = \sum_j f^j[u], \quad \text{where} \quad f^j[u] = c_j(y) \prod_{i=1}^{P_j} \partial_{\vec{x}}^{\vec{p}_i^j} u(\vec{x}, y, t), \tag{7}$$

*for $P_j$ the order of each nonlinear term, and for some $M$-dimensional index $\vec{p}_i^j$. (Sometimes we detail the case when $f[u]$ has only one term in its sum.)*

Roberts (2015) considered systems such as (6) on the multiscale domain $\mathbb{X} \times \mathbb{Y}$ for *one dimensional* $\mathbb{X} \subset \mathbb{R}$, and by analysing the dynamics of the system at some cross-sectional station $X \in \mathbb{X}$, constructed a reduced PDE for the slowly varying dynamics. The construction relied on a Taylor expansion of the field $u$ to order $N$, with the expansion made exact by including the Lagrange's Remainder term in the derivation. Analysis of the Taylor coefficients then reveals the slow behaviour of the system near $X \in \mathbb{X}$ within a centre manifold. Then a projection of the $u$ field PDE onto this centre manifold, generalised to the union over all stations $X \in \mathbb{X}$, defined the slow PDE within domain $\mathbb{X}$, and a projection of Lagrange's Remainder determines the error of the PDE. Roberts and Bunder (2017) analogously considered linear systems on the multiscale domain $\mathbb{X} \times \mathbb{Y}$, but generalised to $M$-dimensional $\mathbb{X} \subset \mathbb{R}^M$. Here we further generalise these earlier developments to the class of nonlinear PDEs (6) that have multiple macroscale dimensions.

## 2.1 Large-scales modulates local Taylor coefficients

In this section we develop the procedure of Roberts and Bunder (2017, §3.1) to the additional complication of nonlinear effects $f[u]$. Both the field $u$ and the nonlinearity $f[u]$ are written as Taylor expansions with Lagrange Remainder terms. These remainder terms ensure the analysis of the dynamics of the system is exact.

The Taylor expansions and subsequent analysis require some assumptions about the smoothness of $u$. For $k_{\max}$ denoting the largest magnitude derivative in the linear term of PDE (6), for $p_i^j$ representing all magnitude derivatives in the nonlinear term (7), and for Taylor expansion to order $N$, the field $u$ must be in differentiability class $C^{N+\max(p_i^j, k_{\max})}$.

At every cross-section station $\vec{X} \in \mathbb{X} \subset \mathbb{R}^M$ we expand the field $u$ as an $N$th order Taylor multinomial about $\vec{x} = \vec{X}$:

$$u(\vec{x}, y, t) = \sum_{|\vec{n}|=0}^{N-1} u^{(\vec{n})}(\vec{X}, y, t) \frac{(\vec{x} - \vec{X})^{\vec{n}}}{\vec{n}!} + \sum_{|\vec{n}|=N} u^{(\vec{n})}(\vec{X}, \vec{x}, y, t) \frac{(\vec{x} - \vec{X})^{\vec{n}}}{\vec{n}!}, \tag{8a}$$

with the multi-index factorial $\vec{n}! := n_1! n_2! \cdots n_M!$, the multi-index magnitude $|\vec{n}| = n_1 + n_2 + \cdots + n_M$, the multi-index power $\vec{x}^{\vec{n}} := x_1^{n_1} x_2^{n_2} \cdots x_M^{n_M}$, and where

- in the first sum, for $|\vec{n}| < N$, the coefficients

$$u^{(\vec{n})}(\vec{X}, y, t) := \partial_{\vec{x}}^{\vec{n}} u \big|_{\vec{x}=\vec{X}}, \qquad (8b)$$

and $u^{(\vec{n})} : \mathbb{X} \times \mathbb{Y} \times \mathbb{T} \to \mathbb{U}$;

- and in the second sum, for $|\vec{n}| = N$, by Lagrange's Remainder Theorem for multivariate Taylor series (Roberts and Bunder 2017, eq. 18, e.g.), the coefficients

$$u^{(\vec{n})}(\vec{X}, \vec{x}, y, t) := N \int_0^1 (1-s)^{N-1} \partial_{\vec{x}}^{\vec{n}} u \big|_{\vec{X}+s(\vec{x}-\vec{X})} \, ds, \qquad (8c)$$

and $u^{(\vec{n})} : \mathbb{X} \times \mathbb{X} \times \mathbb{Y} \times \mathbb{T} \to \mathbb{U}$.

The Taylor expansion of the nonlinear term $f[u]$ in PDE (6) is expressed in the same way as for the field $u$; that is, we do not expand $f[u]$ in a series in $u$, but instead expand $f[u(\vec{x}, y, t)]$ in a series in $\vec{x}$ to order $N$. As Assumption 1 specifies that $f[u]$ is a sum of a product of linear functions, (7), the smoothness requirements for constructing the Taylor expansion of $f[u]$ are satisfied by $u$ being of class $C^{N+\max(p_i^j, k_{\max})}$. The $N$th order Taylor multinomial of $f[u]$ in $\vec{x}$ about $\vec{x} = \vec{X}$ is

$$f[u(\vec{x}, y, t)] = \sum_{|\vec{n}|=0}^{N-1} f^{(\vec{n})}(\vec{X}, y, t) \frac{(\vec{x}-\vec{X})^{\vec{n}}}{\vec{n}!} + \sum_{|\vec{n}|=N} f^{(\vec{n})}(\vec{X}, \vec{x}, y, t) \frac{(\vec{x}-\vec{X})^{\vec{n}}}{\vec{n}!}, \qquad (9a)$$

where

- for $|\vec{n}| < N$,

$$f^{(\vec{n})}(\vec{X}, y, t) := \partial_{\vec{x}}^{\vec{n}} f[u(\vec{x}, y, t)]_{\vec{x}=\vec{X}}; \qquad (9b)$$

- and, for $|\vec{n}| = N$,

$$f^{(\vec{n})}(\vec{X}, \vec{x}, y, t) := N \int_0^1 (1-s)^{N-1} \partial_{\vec{x}}^{\vec{n}} f[u(\vec{x}, y, t)]_{\vec{X}+s(\vec{x}-\vec{X})} \, ds. \qquad (9c)$$

The Taylor coefficients $f^{(\vec{n})} : \mathbb{U} \to \mathbb{U}$ of the nonlinearity are, in principle, functions of the $u$ field Taylor coefficients $u^{(\vec{k})}$ with $|\vec{k}| \leq N$. They may be obtained by substituting the Taylor expansion (8a) of the field $u$ into equations (9b) and (9c). For example, in the case of nonlinearity $f[u]$ being only one term, a direct substitution of (8a) into the nonlinear form (7) gives

$$f[u] = c(y) \prod_{i=1}^{P} \partial_{\vec{x}}^{\vec{p}_i} \left[ \sum_{|\vec{n}|=0}^{N-1} u^{(\vec{n})}(\vec{X}, y, t) \frac{(\vec{x}-\vec{X})^{\vec{n}}}{\vec{n}!} + \sum_{|\vec{n}|=N} u^{(\vec{n})}(\vec{X}, \vec{x}, y, t) \frac{(\vec{x}-\vec{X})^{\vec{n}}}{\vec{n}!} \right]$$

$$= c(y) \prod_{i=1}^{P} \left[ \sum_{|\vec{n}|=0}^{N-1-|\vec{p}_i|} u^{(\vec{n}+\vec{p}_i)}(\vec{X}, y, t) \frac{(\vec{x}-\vec{X})^{(\vec{n}-\vec{p}_i)}}{(\vec{n}-\vec{p}_i)!} \right.$$

$$\left. + \sum_{|\vec{n}|=N-|\vec{p}_i|} \sum_{\vec{m}=\vec{0}}^{\vec{p}_i} \binom{\vec{p}_i}{\vec{m}} \partial_{\vec{x}}^{\vec{m}} u^{(\vec{n}+\vec{p}_i)}(\vec{X}, \vec{x}, y, t) \frac{(\vec{x}-\vec{X})^{(\vec{n}+\vec{m})}}{(\vec{n}+\vec{m})!} \right], \qquad (10)$$

where, for sums over multi-indices such as $\sum_{\vec{m}=\vec{k}}^{\vec{\ell}}$ we require that $k_i \leq m_i \leq \ell_i$ for each component $i = 1, 2, \ldots, M$.

As in the linear case (Roberts and Bunder 2017, Eq. (19)), the multivariate Taylor multinomial (8a) of a field $u$ gives, after some rearrangement, that the $\vec{\ell}$th spatial derivative

$$\partial_{\vec{x}}^{\vec{\ell}} u = \sum_{|\vec{n}|=0}^{N-|\vec{\ell}|-1} u^{(\vec{n}+\vec{\ell})} \frac{(\vec{x}-\vec{X})^{\vec{n}}}{\vec{n}!}$$
$$+ \sum_{|\vec{n}|=N} \sum_{\vec{m}=(\vec{n}-\vec{\ell})^{\oplus}}^{\vec{n}} \binom{\vec{\ell}}{\vec{n}-\vec{m}} \partial_{\vec{x}}^{\vec{m}+\vec{\ell}-\vec{n}} u^{(\vec{n})} \frac{(\vec{x}-\vec{X})^{\vec{m}}}{\vec{m}!}, \quad (11)$$

where appearing for the limits of some sums, $(\vec{k})^{\oplus}$ denotes the multi-index vector with $i$th component $\max(k_i, 0)$, thus ensuring all multi-index components are non-negative in the sums. Now, for every $\vec{\ell}$, $|\vec{\ell}| \leq N$, take the $\vec{\ell}$th spatial derivative of the nonlinear PDE (6),

$$\partial_{\vec{x}}^{\vec{\ell}} \left( \frac{\partial u}{\partial t} \right) = \partial_{\vec{x}}^{\vec{\ell}} \left( \sum_{|\vec{k}|=0}^{\infty} \mathcal{L}_{\vec{k}} \partial_{\vec{x}}^{\vec{k}} u + f[u] \right) \implies \frac{\partial(\partial_{\vec{x}}^{\vec{\ell}} u)}{\partial t} = \sum_{|\vec{k}|=0}^{\infty} \mathcal{L}_{\vec{k}} \partial_{\vec{x}}^{\vec{\ell}+\vec{k}} u + \partial_{\vec{x}}^{\vec{\ell}} f[u],$$

and substitute (11) for the spatial derivatives of field $u$ (replacing $\vec{\ell}$ with $\vec{\ell} + \vec{k}$ where appropriate),

$$\sum_{|\vec{n}|=0}^{N-|\vec{\ell}|-1} \frac{\partial u^{(\vec{n}+\vec{\ell})}}{\partial t} \frac{(\vec{x}-\vec{X})^{\vec{n}}}{\vec{n}!} + \sum_{|\vec{n}|=N} \sum_{\vec{m}=(\vec{n}-\vec{\ell})^{\oplus}}^{\vec{n}} \binom{\vec{\ell}}{\vec{n}-\vec{m}} \partial_{\vec{x}}^{\vec{m}+\vec{\ell}-\vec{n}} \frac{\partial u^{(\vec{n})}}{\partial t} \frac{(\vec{x}-\vec{X})^{\vec{m}}}{\vec{m}!}$$
$$= \sum_{|\vec{k}|=0}^{\infty} \mathcal{L}_{\vec{k}} \sum_{|\vec{n}|=0}^{N-|\vec{\ell}+\vec{k}|-1} u^{(\vec{n}+\vec{\ell}+\vec{k})} \frac{(\vec{x}-\vec{X})^{\vec{n}}}{\vec{n}!}$$
$$+ \sum_{|\vec{k}|=0}^{\infty} \mathcal{L}_{\vec{k}} \sum_{|\vec{n}|=N} \sum_{\vec{m}=(\vec{n}-\vec{\ell}-\vec{k})^{\oplus}}^{\vec{n}} \binom{\vec{\ell}+\vec{k}}{\vec{n}-\vec{m}} \partial_{\vec{x}}^{\vec{m}+\vec{\ell}+\vec{k}-\vec{n}} u^{(\vec{n})} \frac{(\vec{x}-\vec{X})^{\vec{m}}}{\vec{m}!} + \partial_{\vec{x}}^{\vec{\ell}} f[u]. \quad (12)$$

As the multivariate Taylor multinomial (8a) is exact, for all stations $\vec{X} \in \mathbb{X}$ and $\vec{x} \in \mathbb{X}$, equation (12) is exact for every $\vec{x} \in \chi(\vec{X})$, where $\chi(\vec{X})$ is an open subset of $\mathbb{X}$ such that for all points $\vec{x} \in \chi(\vec{X})$ the convex combination $\vec{X} + s(\vec{x}-\vec{X}) \in \chi(\vec{X})$ for every $0 \leq s \leq 1$; this condition ensures that when we take the limit $\vec{x} \to \vec{X}$, $\vec{x}$ will always remain inside $\chi(\vec{X}) \subset \mathbb{X}$ and $(\vec{x}-\vec{X}) \to \vec{0}$.

Now set $\vec{x} = \vec{X}$ in equation (12) so that all terms containing factors of $(\vec{x}-\vec{X})$ vanish. Unless otherwise specified, hereafter $u^{(\vec{n})}$ denotes $u^{(\vec{n})}(\vec{X}, y, t)$ when $|\vec{n}| < N$ and denotes $u^{(\vec{n})}(\vec{X}, \vec{X}, y, t)$ when $|\vec{n}| = N$. Similarly, $f^{(\vec{n})}$ denotes $f^{(\vec{n})}(\vec{X}, y, t)$ when $|\vec{n}| < N$ and denotes $f^{(\vec{n})}(\vec{X}, \vec{X}, y, t)$ when $|\vec{n}| = N$. In addition, swap the $\vec{n}$ and $\vec{\ell}$ multi-indices in (12). Then,

$$\frac{\partial u^{(\vec{n})}}{\partial t} = \sum_{|\vec{k}|=0}^{N-|\vec{n}|} \mathcal{L}_{\vec{k}} u^{(\vec{n}+\vec{k})} + f^{(\vec{n})} + r_{\vec{n}}, \quad \text{for every } |\vec{n}| \leq N, \quad (13a)$$

where the remainder

$$r_{\vec{n}} = \sum_{|\vec{k}| \geq 1} \sum_{\substack{|\vec{\ell}|=N \\ \vec{\ell} \lneq \vec{n}+\vec{k}}} \mathcal{L}_{\vec{k}} \binom{\vec{k}+\vec{n}}{\vec{\ell}} \left[ \partial_{\vec{x}}^{\vec{k}+\vec{n}-\vec{\ell}} u^{(\vec{\ell})}(\vec{X}, \vec{x}, y, t) \right]_{\vec{x}=\vec{X}}, \quad (13b)$$

where throughout $\lneq$ means $\leq$ for each component, but excluding exact equality of the two multi-indices. The second term on the right-hand side of (12) (when $|\vec{k}| \geq 1$) determines the remainder (13b). Since multi-index $\vec{n} \in \mathbb{N}_0^M$ and $|\vec{n}| \leq N$, the total number of unique multi-indices, and thus the number of coupled ODEs (13a), is

$$\mathcal{N} = \binom{N+M}{M}. \tag{14}$$

For all indices $|\vec{n}| = 0, \ldots, N$, the $u^{(\vec{n})}$ terms in equation (13a) are evaluated at station $\vec{X}$, but the spatial derivatives of $u^{(\vec{n})}(\vec{X}, \vec{x}, y, t)$ with $|\vec{n}| = N$ that appear in the remainder term $r_{\vec{n}}$ (13b) couple the dynamics at station $\vec{X}$ to dynamics of the system along the line joining fixed station $\vec{X}$ to variable position $\vec{x}$, that is, the dynamics at $\vec{X}$ are coupled to the dynamics at points in the neighbourhood $\chi(\vec{X})$. This dependence of derivatives of $u^{(\vec{n})}(\vec{X}, \vec{x}, y, t)$ on the dynamics at points in $\chi(\vec{X})$ is directly seen from an application of the integral mean value theorem on equation (8c). By this theorem, there exists some $\hat{s} \in (0, 1)$ such that

$$u^{(\vec{n})}(\vec{X}, \vec{x}, y, t) = N \partial_{\vec{x}}^{\vec{n}} u \big|_{\vec{X} + \hat{s}(\vec{x} - \vec{X})} \int_0^1 (1-s)^{N-1} \, ds = \partial_{\vec{x}}^{\vec{n}} u \big|_{\vec{X} + \hat{s}(\vec{x} - \vec{X})},$$

and $\vec{X} + \hat{s}(\vec{x} - \vec{X}) \in \chi(\vec{X})$. Spatial derivatives of $u^{(\vec{n})}(\vec{X}, \vec{x}, y, t)$ retain dependence on $\hat{s}$, and thus on the dynamics about $\vec{X}$, even when evaluated at $\vec{x} = \vec{X}$. In contrast, $u^{(\vec{n})}(\vec{X}, \vec{X}, y, t) = \partial_{\vec{x}}^{\vec{n}} u \big|_{\vec{X}}$ is independent of $\hat{s}$. Whereas an $\hat{s} \in (0, 1)$ must exist for each $u^{(\vec{n})}(\vec{X}, \vec{x}, y, t)$, these $\hat{s}$ are generally not determined, and so we view gradients of $u^{(\vec{n})}(\vec{X}, \vec{x}, y, t)$ as 'uncertain'. We therefore classify the remainder $r_{\vec{n}}$ as *uncertain forcing* which couple the local dynamics at $\vec{X}$ to the dynamics in its neighbourhood, and thereby to the global dynamics over $\mathbb{X}$.

The nonlinear $f^{(\vec{n})}$ may also contain 'uncertain' gradients of $u^{(\vec{n})}(\vec{X}, \vec{x}, y, t)$, depending on the particular nonlinearity. For example, for the case of a single-term nonlinearity we obtain the last line of equation (10) which contains spatial derivatives up to order $\vec{p}_i$ of $u^{(\vec{n} + \vec{p}_i)}(\vec{X}, \vec{x}, y, t)$ where $|\vec{n} + \vec{p}_i| = N$. So, if at least one $\vec{p}_i > \vec{0}$, the nonlinear term contains uncertain gradients which couple the dynamics at $\vec{X}$ to the dynamics in $\chi(\vec{X})$. Section 2.2 explicitly identifies these uncertain gradients in the nonlinear $f^{(\vec{n})}$.

## 2.2 Generating multinomial and PDE

We now pack all the multivariate Taylor coefficients $u^{(\vec{n})}$ together into a generating function (multinomial). For every station $\vec{X} \in \mathbb{X}$ and time $t \in \mathbb{T}$ consider the field $u$ in terms of a local Taylor multinomial (8a) about the cross-section $\vec{x} = \vec{X}$. In terms of the indeterminate $\vec{\xi} \in \mathbb{R}^M$, define the generating multinomial

$$\tilde{u}(\vec{X}, t) := \sum_{|\vec{n}|=0}^{N-1} \frac{\vec{\xi}^{\vec{n}}}{\vec{n}!} u^{(\vec{n})}(\vec{X}, y, t) + \sum_{|\vec{n}|=N} \frac{\vec{\xi}^{\vec{n}}}{\vec{n}!} u^{(\vec{n})}(\vec{X}, \vec{X}, y, t), \tag{15}$$

where this generating multinomial $\tilde{u}$, through its range denoted by $\mathbb{U}_N$, is implicitly a function of the indeterminate $\vec{\xi}$ and the cross-sectional variable $y$. This generating multinomial $\tilde{u} : \mathbb{X} \times \mathbb{T} \to \mathbb{U}_N$ for the vector space

$$\mathbb{U}_N := \mathbb{U} \otimes_t \mathbb{G}_N \quad \text{where} \quad \mathbb{G}_N := \{\text{multinomials in } \vec{\xi} \text{ of degree} \leq N\},$$

and where $\otimes_t$ represents the vector space tensor product. The generating operator

$$\mathcal{G} = \left[ \sum_{|\vec{n}|=0}^N \frac{\vec{\xi}^{\vec{n}}}{\vec{n}!} \partial_{\vec{x}}^{\vec{n}} \right]_{\vec{x}=\vec{X}}, \tag{16}$$

acts to convert the original field $u(\vec{x}, y, t)$ into the generating multinomial $\tilde{u}(\vec{X}, t) = \mathcal{G}u(\vec{x}, y, t)$. The generating operator (16) similarly converts the nonlinear term of PDE (6) into a multinomial in $\vec{\xi}$,

$$\tilde{f}[\tilde{u}] := \mathcal{G}f[u(\vec{x}, y, t)] = \sum_{|\vec{n}|=0}^{N-1} \frac{\vec{\xi}^{\vec{n}}}{\vec{n}!} f^{(\vec{n})}(\vec{X}, y, t) + \sum_{|\vec{n}|=N} \frac{\vec{\xi}^{\vec{n}}}{\vec{n}!} f^{(\vec{n})}(\vec{X}, \vec{X}, y, t), \quad (17)$$

with $\tilde{f}[\tilde{u}] : \mathbb{U}_N \to \mathbb{U}_N$ appearing as the nonlinear term in (18) of the next Proposition 2.

We introduce the multinomial $\tilde{u}$ because it is more conveniently compact to deal with one multinomial and one PDE than the $N$ Taylor coefficients $u^{(\vec{n})}$ and the $N$ differential equations (13) derived in the previous section. Roberts and Bunder (2017, §3.2) constructed a similar multinomial $\tilde{u}$ and PDE via the generating operator $\mathcal{G}$; however, here we make new special provisions for the nonlinear term $f[u]$. Although the compact form of multinomial $\tilde{u}$ is useful, a more important property is the equivalence of the dynamics of $\tilde{u}$ and the original field $u(\vec{x}, y, t)$, as described by Proposition 2.

**Proposition 2.** *Let $u(\vec{x}, y, t)$ be governed by the specified nonlinear PDE (6). Then the dynamics at every locale $\vec{X} \in \mathbb{X} \subset \mathbb{R}^M$ is equivalently governed by the nonlinear PDE*

$$\frac{\partial \tilde{u}}{\partial t} = \tilde{\mathcal{L}}\tilde{u} + \tilde{f}[\tilde{u}] + \tilde{r}[u], \quad (18)$$

*for the generating function multinomial $\tilde{u}(\vec{X}, y, t)$ defined in (15), the 'uncertain' forcing $\tilde{r}[u]$ given by (22), nonlinear $\tilde{f}[\tilde{u}]$ defined by (17), and operator*

$$\tilde{\mathcal{L}} = \sum_{|\vec{k}|=0}^{N} \mathfrak{L}_{\vec{k}} \partial_{\vec{\xi}}^{\vec{k}}. \quad (19)$$

To establish Proposition 2, we first show that multinomial $\tilde{u}(\vec{X}, y, t)$ (15) satisfies PDE (18). To construct a PDE for $\tilde{u}$, take the time derivative of (15) and replace $\partial u^{(\vec{n})}/\partial t$ using (13a):

$$\frac{\partial \tilde{u}}{\partial t} = \sum_{|\vec{n}|=0}^{N} \frac{\vec{\xi}^{\vec{n}}}{\vec{n}!} \left[ \sum_{|\vec{k}|=0}^{N-|\vec{n}|} \mathfrak{L}_{\vec{k}} u^{(\vec{n}+\vec{k})} \right] + \sum_{|\vec{n}|=0}^{N} \frac{\vec{\xi}^{\vec{n}}}{\vec{n}!} f^{(\vec{n})} + \sum_{|\vec{n}|=0}^{N} \frac{\vec{\xi}^{\vec{n}}}{\vec{n}!} r_{\vec{n}}$$

$$= \sum_{|\vec{k}|=0}^{N} \mathfrak{L}_{\vec{k}} \partial_{\vec{\xi}}^{\vec{k}} \tilde{u} + \tilde{f}[\tilde{u}] + \sum_{|\vec{n}|=0}^{N} \frac{\vec{\xi}^{\vec{n}}}{\vec{n}!} r_{\vec{n}}, \quad (20)$$

where in the first term the $\vec{n}$ and $\vec{k}$ sums are exchanged, and we then simplify this term using the useful identity

$$\partial_{\vec{\xi}}^{\vec{k}} \tilde{u} = \sum_{|\vec{n}|=0}^{N-|\vec{k}|} \frac{\vec{\xi}^{\vec{n}}}{\vec{n}!} u^{(\vec{n}+\vec{k})}, \quad (21)$$

obtained from derivatives of the generating multinomial (15) with respect to $\vec{\xi}$. The above PDE (20) is precisely PDE (18) of Proposition 2 with forcing 'remainder' term

$$\tilde{r}[u] = \sum_{|\vec{n}|=0}^{N} \frac{\vec{\xi}^{\vec{n}}}{\vec{n}!} r_{\vec{n}}$$

$$= \sum_{|\vec{k}|\geqslant 1} \mathcal{L}_{\vec{k}} \sum_{|\vec{n}|=0}^{N} \frac{\vec{\xi}^{\vec{n}}}{\vec{n}!} \sum_{\substack{|\vec{\ell}|=N \\ \vec{\ell} \lneq \vec{n}+\vec{k}}} \binom{\vec{k}+\vec{n}}{\vec{\ell}} \left[\partial_{\vec{x}}^{\vec{k}+\vec{n}-\vec{\ell}} u^{(\vec{\ell})}(\vec{X},\vec{x},y,t)\right]_{\vec{x}=\vec{X}}, \qquad (22)$$

upon using expression (13b) for $r_{\vec{n}}$.

The second task for establishing Proposition 2 is to show that the generating PDE (18) and the original PDE (6) describe the same dynamics at every locale $\vec{X} \in \mathbb{X} \subset \mathbb{R}^M$. We do this by providing a more physical interpretation of the generating operator $\mathcal{G}$ and the generating multinomial $\tilde{u}(\vec{X},t)$, beyond just a convenient way to pack the Taylor coefficients of $u(\vec{x},y,t)$.

Consider the Taylor expansion of some general function $g(\vec{x}) \in C^{N+1}$ at $\vec{x} = \vec{X} + \vec{\xi}$ about $\vec{x} = \vec{X}$:

$$[g(\vec{x})]_{\vec{x}=\vec{X}+\vec{\xi}} = g(\vec{X}+\vec{\xi}) = \sum_{|\vec{n}|=0}^{N} \frac{\vec{\xi}^{\vec{n}}}{\vec{n}!} \partial_{\vec{x}}^{\vec{n}} g(\vec{X}) + R_N(g) = \mathcal{G}g(\vec{x}) + \mathcal{O}(|\vec{\xi}|^{N+1}), \quad (23)$$

where $R_N(g)$ is the order $N$ Lagrange remainder term of $g(\vec{X}+\vec{\xi})$ (Roberts and Bunder 2017). So, $\mathcal{G}g(\vec{x})$ evaluates $g(\vec{X}+\vec{\xi})$ correct to $\mathcal{O}(|\vec{\xi}|^{N+1})$. Similarly, $\tilde{u}(\vec{X},t) = \mathcal{G}u(\vec{x},y,t)$ evaluates $u(\vec{X}+\vec{\xi},y,t)$ correct to $\mathcal{O}(|\vec{\xi}|^{N+1})$. We interpret $\tilde{u}(\vec{X},t)$ as the projection of $u(\vec{x},y,t)$ at $\vec{x}=\vec{X}+\vec{\xi}$ onto the space $\mathbb{U}_N = \mathbb{U} \otimes_t \mathbb{G}_N$, with $\mathcal{O}(|\vec{\xi}|^{N+1})$ interpreted not as an error but as the difference between $u(\vec{X}+\vec{\xi},y,t)$ and its projection onto $\mathbb{U}_N$ (Roberts and Bunder 2017). As $\mathcal{G}$ commutes with the temporal derivative $\mathcal{G}\partial u(\vec{x},y,t)/\partial t = \partial \tilde{u}(\vec{X},t)/\partial t$ and $\partial \tilde{u}(\vec{X},t)/\partial t$ is equivalent to the Taylor expansion of $\partial u(\vec{x},y,t)/\partial t$ at $\vec{x}=\vec{X}+\vec{\xi}$ correct to $\mathcal{O}(|\vec{\xi}|^{N+1})$. Therefore the generating PDE (18) for multinomial $\tilde{u}(\vec{X},t)$ is equivalent to the PDE (6) for $u(\vec{x},y,t)$ evaluated at $\vec{x}=\vec{X}+\vec{\xi}$ correct to $\mathcal{O}(|\vec{\xi}|^{N+1})$. Thus the dynamics of PDE (18) are identical to the dynamics of PDE (6) at every $\vec{x}=\vec{X} \in \mathbb{X}$. This completes the proof of Proposition 2.

The dynamics of the original nonlinear PDE (6) for field $u$ are equivalent to the dynamics of the nonlinear PDE (18) for the $N$ dimensional multinomial $\tilde{u}$ (15); furthermore, the two PDEs are symbolically the same with $u \leftrightarrow \tilde{u}$ and $\vec{x} \leftrightarrow \vec{\xi}$, plus a forcing term. But the advantage of the multinomial form is that the derivatives $\partial_{\vec{\xi}}$ operate only on $\mathbb{G}_N$, that is, multinomials of at most degree $N$ in $\vec{\xi} \in \mathbb{R}^M$, and are thus bounded in $\mathbb{G}_N$. In contrast, the derivatives $\partial_{\vec{x}}$ in the original PDE are potentially unbounded (e.g., for $u$ rapidly oscillating or containing irrational functions). The slowly varying modelling of Section 3.2 takes advantage of the near symbolic equivalence between PDE (6) and PDE (18) with $u \leftrightarrow \tilde{u}$ and $\vec{x} \leftrightarrow \vec{\xi}$.

We now expand the nonlinear term (17) of PDE (18) explicitly in terms of generating multinomial $\tilde{u}$ and nonlinear 'uncertain' terms involving gradients of $u^{(\vec{n})}$ with $|\vec{n}| > N$. Section 3.2 makes use of this expansion to simplify the remainder term $\rho$ of the slow PDE, and the appendix applies the expansion in the construct the slow PDE for the fluid flow example discussed in Section 1.1.

The nonlinear term (17) in the generating PDE (18), expanded according to (7) in Assumption 1 is

$$\tilde{f}[\tilde{u}] = \sum_{j} \sum_{|\vec{n}|=0}^{N} \frac{\vec{\xi}^{\vec{n}}}{\vec{n}!} \left[\partial_{\vec{x}}^{\vec{n}} c_j(y) \prod_{i=1}^{P_j} \partial_{\vec{x}}^{\vec{p}_i^j} u(\vec{x},y,t)\right]_{\vec{x}=\vec{X}}$$

$$= \sum_{j} \sum_{|\vec{n}|=0}^{N} \sum_{\sum_{i=1}^{P_j} \vec{m}_i = \vec{n}} c_j(y) \prod_{i=1}^{P_j} \frac{\vec{\xi}^{\vec{m}_i}}{\vec{m}_i!} \left[\partial_{\vec{x}}^{\vec{m}_i} \partial_{\vec{x}}^{\vec{p}_i^j} u(\vec{x},y,t)\right]_{\vec{x}=\vec{X}}$$

$$
\begin{aligned}
= &\sum_j \sum_{|\vec{n}|=0}^{N} \sum_{\substack{\sum_{i=1}^{P_j} \vec{m}_i = \vec{n} \\ |\vec{m}_i + \vec{p}_i^{\,j}| \leqslant N}} c_j(y) \prod_{i=1}^{P_j} \frac{\vec{\xi}^{\vec{m}_i}}{\vec{m}_i!} \left[ \partial_{\vec{x}}^{\vec{m}_i + \vec{p}_i^{\,j}} u(\vec{x}, y, t) \right]_{\vec{x} = \vec{X}} \\
&+ \sum_j \sum_{|\vec{n}|=0}^{N} \sum_{\substack{\sum_{i=1}^{P_j} \vec{m}_i = \vec{n} \\ \exists |\vec{m}_i + \vec{p}_i^{\,j}| > N}} c_j(y) \prod_{i=1}^{P_j} \frac{\vec{\xi}^{\vec{m}_i}}{\vec{m}_i!} \left[ \partial_{\vec{x}}^{\vec{m}_i + \vec{p}_i^{\,j}} u(\vec{x}, y, t) \right]_{\vec{x} = \vec{X}}. \quad (24)
\end{aligned}
$$

The components with $|\vec{m}_i + \vec{p}_i^{\,j}| > N$ in the last term on the right hand side are 'uncertain', similar to the uncertain forcing $\tilde{r}[u]$ (22), although in the special case where all $\vec{p}_i^{\,j} = \vec{0}$, no such uncertain nonlinear terms exist. For the uncertain gradients, consider expansion (11) with $|\vec{\ell}| > N$ evaluated at $\vec{x} = \vec{X}$,

$$
\partial_{\vec{x}}^{\vec{\ell}} [u(\vec{x}, y, t)]_{\vec{x} = \vec{X}} = \sum_{|\vec{n}| = N,\, \vec{n} \lneq \vec{\ell}} \binom{\vec{\ell}}{\vec{n}} \left[ \partial_{\vec{x}}^{\vec{\ell} - \vec{n}} u^{(\vec{n})}(\vec{X}, \vec{x}, y, t) \right]_{\vec{x} = \vec{X}}.
$$

Using this expansion for the uncertain terms, as well as (8b) and (21) and Assumption 1, we rewrite the nonlinear term (24) as

$$
\begin{aligned}
\tilde{f}[\tilde{u}] = &\sum_j \sum_{|\vec{n}|=0}^{N} \sum_{\substack{\sum_{i=1}^{P} \vec{m}_i = \vec{n} \\ |\vec{m}_i + \vec{p}_i^{\,j}| \leqslant N}} c_j(y) \prod_{i=1}^{P_j} \frac{\vec{\xi}^{\vec{m}_i}}{\vec{m}_i!} \left[ \partial_{\vec{\xi}}^{\vec{m}_i + \vec{p}_i^{\,j}} \tilde{u} \right]_{\vec{\xi} = \vec{0}} \\
&+ \sum_j \sum_{|\vec{n}|=0}^{N} \sum_{\substack{\sum_{i=1}^{P_j} \vec{m}_i = \vec{n} \\ \exists |\vec{m}_i + \vec{p}_i^{\,j}| > N}} c_j(y) \prod_{i=1}^{P_j} \frac{\vec{\xi}^{\vec{m}_i}}{\vec{m}_i!} f_i^j[u, \tilde{u}], \quad (25)
\end{aligned}
$$

where in the second term $f_i[u, \tilde{u}]$ is either a function of the generating multinomial $\tilde{u}$ or of the uncertain gradients of original field $u$,

$$
f_i^j[u, \tilde{u}] := \begin{cases} \left[ \partial_{\vec{\xi}}^{\vec{m}_i + \vec{p}_i^{\,j}} \tilde{u} \right]_{\vec{\xi} = \vec{0}} & \text{for } |\vec{m}_i + \vec{p}_i^{\,j}| \leqslant N, \\ \displaystyle\sum_{\substack{|\vec{k}| = N \\ \vec{k} \lneq \vec{m}_i + \vec{p}_i^{\,j}}} \binom{\vec{m}_i + \vec{p}_i^{\,j}}{\vec{k}} \left[ \partial_{\vec{x}}^{\vec{m}_i + \vec{p}_i^{\,j} - \vec{k}} u^{(\vec{k})}(\vec{X}, \vec{x}, y, t) \right]_{\vec{x} = \vec{X}} & \text{for } |\vec{m}_i + \vec{p}_i^{\,j}| > N. \end{cases}
$$

Thus in the Taylor expansion (17),

$$
\begin{aligned}
f^{(\vec{n})} = &\vec{n}! \sum_j \sum_{\substack{\sum_{i=1}^{P_j} \vec{m}_i = \vec{n} \\ |\vec{m}_i + \vec{p}_i^{\,j}| \leqslant N}} c_j(y) \prod_{i=1}^{P_j} \frac{1}{\vec{m}_i!} \left[ \partial_{\vec{\xi}}^{\vec{m}_i + \vec{p}_i^{\,j}} \tilde{u} \right]_{\vec{\xi} = \vec{0}} \\
&+ \vec{n}! \sum_j \sum_{\substack{\sum_{i=1}^{P_j} \vec{m}_i = \vec{n} \\ \exists |\vec{m}_i + \vec{p}_i^{\,j}| > N}} c_j(y) \prod_{i=1}^{P_j} \frac{1}{\vec{m}_i!} f_i^j[u, \tilde{u}], \quad (26)
\end{aligned}
$$

where the second term contains all uncertain gradients.

# 3 A slow nonlinear model emerges

Section 3.1 constructs the eigenspace which describes the emergent slow dynamics of the generating PDE (18) by analysing a linearisation of the PDE (18). We then show how this eigenspace and associated eigenvalues describe the slow dynamics of original nonlinear PDE (6). It is possible to determine the slow dynamics of PDE (6) from the eigenspace of the linearised Taylor coefficient PDEs (13a), without introducing the generating multinomial and generating function, but then one must explicitly deal with $N$ Taylor coefficients and their coupled $N$ PDEs, as seen in the linear example of Roberts and Bunder (2017) [§2.2]. Employing the linear eigenspace as a foundation to describe the dynamics of a nonlinear system is justified by centre manifold theory (Carr 1981; Aulbach and Wanner 2000; Haragus and Iooss 2011, e.g.) which assures us that generically the stability properties of a linear system with centre-stable dynamics persist under nonlinear perturbations and time-dependent forcing. Section 3.2 extends the analysis to the nonlinear PDE (18) to describe the slow dynamics of the nonlinear system, including coupling between the centre and stable subspaces via uncertain terms in both the forcing and the nonlinear term.

The slow dynamics are characterised by a set of generalised eigenvectors $\tilde{\mathcal{V}}^{\vec{n}}$ which span the centre subspace on which the slow dynamics of $\tilde{u}$ evolve. As these generalised eigenvectors are determined from the linearised PDE, they are the same as those we determined (Roberts and Bunder 2017, §3.3) for linear PDEs. However, there we constructed the generalised eigenvectors $\tilde{\mathcal{V}}^{\vec{n}}$ via a two-step process, firstly constructing generalised eigenvectors for the linearised version of the original PDE (6) (i.e., for $f[u] = 0$), and then mapping these eigenvectors into $\mathbb{U}_N$ (Roberts and Bunder 2017, eq. (37) and Lemma 6). Here we show how to construct the generalised eigenvectors in $\mathbb{U}_N$ directly from the generating PDE (18).

To analyse the eigenspace and determine the slow dynamics of the linearised generating PDE (18), we apply Assumption 3 which describes the eigenspace of $\mathfrak{L}_{\vec{0}}$, the lowest order operator in $\tilde{\mathcal{L}}$. However, Assumption 3 does not provide *necessary* assumptions for the extraction of a slow model; for example, here we derive the slow model after assuming the Hilbert space $\mathbb{U}$ is a centre-stable subspace, but an analogous derivation is possible when $\mathbb{U}$ is a slow-stable subspace (the shallow fluid example of Section 1.1 is on a slow-stable subspace and the code in Appendices A.2 and A.3 permit either slow-stable or centre-stable dynamics).

**Assumption 3.** *We assume the following for the primary case of purely centre-stable dynamics.*

1. *The Hilbert space $\mathbb{U}$ is the direct sum of two closed $\mathfrak{L}_{\vec{0}}$-invariant subspaces, $\mathbb{E}_c^0$ and $\mathbb{E}_s^0$, and the corresponding restrictions of $\mathfrak{L}_{\vec{0}}$ generate strongly continuous semigroups (Gallay 1993; Aulbach and Wanner 1996).*

2. *The operator $\mathfrak{L}_{\vec{0}}$ has a discrete spectrum of eigenvalues $\lambda_1, \lambda_2, \ldots$ (repeated according to multiplicity) with corresponding linearly independent (possibly generalised) eigenvectors $v_1^{\vec{0}}, v_2^{\vec{0}}, \ldots$ that are complete ($\mathbb{U} = \text{span}\{v_1^{\vec{0}}, v_2^{\vec{0}}, \ldots\}$).*

3. *The first $\mathfrak{m}$ eigenvalues $\lambda_1, \ldots, \lambda_{\mathfrak{m}}$ of $\mathfrak{L}_{\vec{0}}$ all have real part satisfying $|\Re\lambda_j| \leq \alpha$ and hence the $\mathfrak{m}$-dimensional centre subspace $\mathbb{E}_c^0 = \text{span}\{v_1^{\vec{0}}, \ldots, v_{\mathfrak{m}}^{\vec{0}}\}$ (Chicone 2006, Chap. 4, e.g.).*

4. *All other eigenvalues $\lambda_{\mathfrak{m}+1}, \lambda_{\mathfrak{m}+2}, \ldots$ have real part negative and well separated from the centre eigenvalues, namely $\Re\lambda_j \leq -\beta < -N\alpha$ for $j =$* $\mathfrak{m}+1, \mathfrak{m}+2, \ldots$, and so the stable subspace $\mathbb{E}_s^0 = \text{span}\{v_{\mathfrak{m}+1}^{\vec{0}}, v_{\mathfrak{m}+2}^{\vec{0}}, \ldots\}$. For clarity, say the number of stable eigenvalues is $\mathfrak{m}'$, so that the stable subspace $\mathbb{E}_s^0$ is $\mathfrak{m}'$-dimensional, although the number of stable eigenvalues may be infinite, $\mathfrak{m}' \to \infty$.

For convenience, Definition 4 packs the $\mathfrak{m}$ eigenvectors which span the centre subspace $\mathbb{E}_c^0$ of $\mathfrak{L}_{\vec{0}}$ into one matrix $V^{\vec{0}}$, and similarly packs the eigenvalues into the matrix $A_{\vec{0}}$.

**Definition 4.** *Assumption 3 identifies a subset of $\mathfrak{m}$ eigenvectors of $\mathfrak{L}_{\vec{0}}$ which span the centre subspace $\mathbb{E}_c^0 \subset \mathbb{U}$.*

- *With these eigenvectors define*

$$V^{\vec{0}} := \begin{bmatrix} v_1^{\vec{0}} & v_2^{\vec{0}} & \cdots & v_\mathfrak{m}^{\vec{0}} \end{bmatrix} \in \mathbb{U}^{1 \times \mathfrak{m}}.$$

- *Since the centre subspace is an invariant space of $\mathfrak{L}_{\vec{0}}$, define complex matrix $A_{\vec{0}} \in \mathbb{C}^{\mathfrak{m} \times \mathfrak{m}}$ to be such that $\mathfrak{L}_{\vec{0}} V^{\vec{0}} = V^{\vec{0}} A_{\vec{0}}$ (often $A_{\vec{0}}$ will be in Jordan form, but it is not necessarily so).*

- *Use $\langle \cdot, \cdot \rangle$ to also denote the inner product on the Hilbert space $\mathbb{U}$, $\langle \cdot, \cdot \rangle : \mathbb{U} \times \mathbb{U} \to \mathbb{C}$, the field of complex numbers.*

  *Interpret this inner product when acting on two matrices/vectors with elements in $\mathbb{U}$ as the matrix/vector of the corresponding elementwise inner products. For example, for $Z^{\vec{0}}, V^{\vec{0}} \in \mathbb{U}^{1 \times \mathfrak{m}}$, $\langle Z^{\vec{0}}, V^{\vec{0}} \rangle \in \mathbb{C}^{\mathfrak{m} \times \mathfrak{m}}$.*

- *Define $Z^{\vec{0}} \in \mathbb{U}^{1 \times \mathfrak{m}}$ to have $\mathfrak{m}$ linearly independent columns which are the $\mathfrak{m}$ left eigenvectors of $\mathfrak{L}_{\vec{0}}$, ordered such that the $j$th columns of $V^{\vec{0}}$ and $Z^{\vec{0}}$ have the same eigenvalue and normalised such that $\langle Z^{\vec{0}}, V^{\vec{0}} \rangle = I_\mathfrak{m}$.*

The next Section 3.1 uses the centre subspace eigenvectors $V^{\vec{0}}$ of $\mathfrak{L}_{\vec{0}}$ to generate a set of generalised eigenvectors of $\tilde{\mathcal{L}}$ which describe the slow dynamics of linear PDE $\partial_t \tilde{u} = \tilde{\mathcal{L}} \tilde{u}$ confined to the centre subspace.

## 3.1 Generalised eigenvectors span the centre subspace

We invoke Assumption 3 to construct a set of eigenvectors (possibly generalised) which span the centre subspace $\mathbb{E}_c^N \subset \mathbb{U}_N$ of the linear operator $\tilde{\mathcal{L}}$ (19). These eigenvectors capture the slow behaviour of the linear PDE

$$\frac{\partial \tilde{u}}{\partial t} = \tilde{\mathcal{L}} \tilde{u}, \tag{27}$$

which is the linearisation of generating PDE (18), with neglected forcing terms.

For $0 < |\vec{n}| \leqslant N$, we construct the generalised eigenvector $\tilde{\mathcal{V}}^{\vec{n}} \in \mathbb{U}^{1 \times \mathfrak{m}} \otimes_t \mathbb{G}_N = \mathbb{U}_N^{1 \times \mathfrak{m}}$ from the following recurrence relations, beginning with $\tilde{\mathcal{V}}^{\vec{0}} = V^{\vec{0}}$,

$$A_{\vec{n}} := \sum_{0 < |\vec{k}|, \vec{k} \leqslant \vec{n}} \langle Z^{\vec{0}}, \mathfrak{L}_{\vec{k}} \tilde{\mathcal{V}}^{\vec{n}-\vec{k}} \rangle_{\vec{\xi}=\vec{0}}, \tag{28a}$$

$$\mathfrak{L}_{\vec{0}} \tilde{\mathcal{V}}^{\vec{n}} - \tilde{\mathcal{V}}^{\vec{n}} A_{\vec{0}} = - \sum_{0 < |\vec{k}|, \vec{k} \leqslant \vec{n}} \mathfrak{L}_{\vec{k}} \tilde{\mathcal{V}}^{\vec{n}-\vec{k}} + \sum_{0 < |\vec{k}|, \vec{k} \leqslant \vec{n}} \tilde{\mathcal{V}}^{\vec{n}-\vec{k}} A_{\vec{k}}, \tag{28b}$$

$$\langle Z^{\vec{0}}, \tilde{\mathcal{V}}^{\vec{n}} \rangle = \frac{\vec{\xi}^{\vec{n}}}{\vec{n}!} I_\mathfrak{m}. \tag{28c}$$

The $\mathfrak{m}$ rows of all $\tilde{\mathcal{V}}^{\vec{\mathfrak{n}}}$ with $|\vec{\mathfrak{n}}| \leqslant \mathcal{N}$ form a subset of $\mathbb{U}_\mathcal{N}$ with $\mathfrak{m}\mathcal{N}$ elements. Below we show that these $\mathfrak{m}\mathcal{N}$ elements are generalised eigenvectors of $\tilde{\mathcal{L}}$ which span the centre subspace $\mathbb{E}_c^\mathcal{N}$. To do this we show that the $\mathfrak{m}\mathcal{N}$ elements are linearly independent and that the generalised eigenvector equation $\tilde{\mathcal{L}}\tilde{\mathcal{V}}^{\vec{\mathfrak{n}}} - \tilde{\mathcal{V}}^{\vec{\mathfrak{n}}} A_{\vec{0}}$ only produces linear combinations of $\tilde{\mathcal{V}}^{\vec{\mathfrak{k}}}$ with $\vec{0} \leqslant \vec{\mathfrak{k}} < \vec{\mathfrak{n}}$.

The inner product (28c) ensures that $\tilde{\mathcal{V}}^{\vec{\mathfrak{n}}} = \frac{\vec{\xi}^{\vec{\mathfrak{n}}}}{\vec{\mathfrak{n}}!}V^{\vec{0}} + \tilde{V}^{\vec{\mathfrak{n}}}$ for some $\tilde{V}^{\vec{\mathfrak{n}}} \in \mathbb{U}_\mathcal{N}^{1 \times \mathfrak{m}}$ such that $\langle Z^{\vec{0}}, \tilde{V}^{\vec{\mathfrak{n}}} \rangle = 0_\mathfrak{m}$ for all $|\vec{\mathfrak{n}}| > 0$. Further, since the $\frac{\vec{\xi}^{\vec{\mathfrak{n}}}}{\vec{\mathfrak{n}}!}V^{\vec{0}}$ part of $\tilde{\mathcal{V}}^{\vec{\mathfrak{n}}}$ gives zero in the left hand side of (28b), the objective of (28b) is to determine the $\tilde{V}^{\vec{\mathfrak{n}}}$ part of $\tilde{\mathcal{V}}^{\vec{\mathfrak{n}}}$. As the right hand side of (28b) is a function of $\tilde{\mathcal{V}}^{\vec{\mathfrak{k}}}$ with $\vec{0} \leqslant \vec{\mathfrak{k}} < \vec{\mathfrak{n}}$, we conclude that $\tilde{V}^{\vec{\mathfrak{n}}}$ has order of $\vec{\xi}$ no larger than the order of these $\tilde{\mathcal{V}}^{\vec{\mathfrak{k}}}$. Since we know that $V^{\vec{0}}$ is independent of $\vec{\xi}$, for $|\vec{\mathfrak{n}}| = 1$ equation (28b) ensures that $\tilde{V}^{\vec{\mathfrak{n}}}$ is independent of $\vec{\xi}$ and the highest order of $\tilde{\mathcal{V}}^{\vec{\mathfrak{n}}}$ when $|\vec{\mathfrak{n}}| = 1$ must be $\vec{\xi}^{\vec{\mathfrak{n}}}$. By induction we conclude that for any $\vec{\mathfrak{n}}$, $\tilde{V}^{\vec{\mathfrak{n}}}$ is of order $\vec{\mathfrak{k}}$ in $\vec{\xi}$, where $\vec{0} \leqslant \vec{\mathfrak{k}} < \vec{\mathfrak{n}}$, and $\tilde{\mathcal{V}}^{\vec{\mathfrak{n}}}$ is of order $\vec{\mathfrak{n}}$ in $\vec{\xi}$. Thus $\tilde{\mathcal{V}}^{\vec{\mathfrak{n}}}$ is an $\vec{\mathfrak{n}}$th order multinomial in $\mathbb{U}_\mathcal{N}^{1 \times \mathfrak{m}}$ and for all $|\vec{\mathfrak{n}}| \leqslant \mathcal{N}$ we have $\mathcal{N}$ linearly independent $\tilde{\mathcal{V}}^{\vec{\mathfrak{n}}}$.

Now consider the rows of each $\tilde{\mathcal{V}}^{\vec{\mathfrak{n}}}$. Since

$$\tilde{\mathcal{V}}^{\vec{\mathfrak{n}}} = \frac{\vec{\xi}^{\vec{\mathfrak{n}}}}{\vec{\mathfrak{n}}!}V^{\vec{0}} + \tilde{V}^{\vec{\mathfrak{n}}} = \frac{\vec{\xi}^{\vec{\mathfrak{n}}}}{\vec{\mathfrak{n}}!}\begin{bmatrix} v_1^{\vec{0}} & v_2^{\vec{0}} & \cdots & v_\mathfrak{m}^{\vec{0}} \end{bmatrix} + \tilde{V}^{\vec{\mathfrak{n}}},$$

with linearly independent eigenvectors $v_j^{\vec{0}}$ for $j = 1, \ldots, \mathfrak{m}$, and since $\langle Z^{\vec{0}}, \tilde{V}^{\vec{\mathfrak{n}}} \rangle = 0_\mathfrak{m}$, each of the $\mathfrak{m}$ elements of $\tilde{\mathcal{V}}^{\vec{\mathfrak{n}}}$ are linearly independent. Therefore, the $\mathfrak{m}$ elements of all $\tilde{\mathcal{V}}^{\vec{\mathfrak{n}}}$ with $|\vec{\mathfrak{n}}| \leqslant \mathcal{N}$, form a set of $\mathfrak{m}\mathcal{N}$ linearly independent elements of $\mathbb{U}_\mathcal{N}$.

To show that the rows of $\tilde{\mathcal{V}}^{\vec{\mathfrak{n}}}$ are generalised eigenvectors of $\tilde{\mathcal{L}}$ in the centre subspace $\mathbb{E}_c^\mathcal{N}$, consider

$$\begin{aligned}
\tilde{\mathcal{L}}\tilde{\mathcal{V}}^{\vec{\mathfrak{n}}} - \tilde{\mathcal{V}}^{\vec{\mathfrak{n}}} A_{\vec{0}} &= \sum_{|\vec{\mathfrak{k}}|=0}^{\mathcal{N}} \mathfrak{L}_{\vec{\mathfrak{k}}} \partial_{\vec{\xi}}^{\vec{\mathfrak{k}}} \tilde{\mathcal{V}}^{\vec{\mathfrak{n}}} - \tilde{\mathcal{V}}^{\vec{\mathfrak{n}}} A_{\vec{0}} \\
&= \sum_{0 \leqslant \vec{\mathfrak{k}} \leqslant \vec{\mathfrak{n}}} \mathfrak{L}_{\vec{\mathfrak{k}}} \partial_{\vec{\xi}}^{\vec{\mathfrak{k}}} \tilde{\mathcal{V}}^{\vec{\mathfrak{n}}} - \tilde{\mathcal{V}}^{\vec{\mathfrak{n}}} A_{\vec{0}} \quad \text{since } \tilde{\mathcal{V}}^{\vec{\mathfrak{n}}} \text{ is order } \vec{\mathfrak{n}} \text{ in } \vec{\xi} \\
&= \sum_{0 \leqslant \vec{\mathfrak{k}} \leqslant \vec{\mathfrak{n}}} \mathfrak{L}_{\vec{\mathfrak{k}}} \tilde{\mathcal{V}}^{\vec{\mathfrak{n}}-\vec{\mathfrak{k}}} - \tilde{\mathcal{V}}^{\vec{\mathfrak{n}}} A_{\vec{0}} \quad \text{from Lemma 5} \\
&= \sum_{0 < \vec{\mathfrak{k}} \leqslant \vec{\mathfrak{n}}} \tilde{\mathcal{V}}^{\vec{\mathfrak{n}}-\vec{\mathfrak{k}}} A_{\vec{\mathfrak{k}}} \quad \text{from rearranging (28b)}. \quad (29)
\end{aligned}$$

The left-hand side only produces $\tilde{\mathcal{V}}^{\vec{\mathfrak{k}}}$ with $\vec{0} \leqslant \vec{\mathfrak{k}} < \vec{\mathfrak{n}}$, and thus the rows of $\tilde{\mathcal{V}}^{\vec{\mathfrak{n}}}$ are generalised eigenvectors of rank $\vec{\mathfrak{n}}$ with eigenvalues in matrix $A_{\vec{0}}$. Since the rows of all $\tilde{\mathcal{V}}^{\vec{\mathfrak{n}}}$ with $|\vec{\mathfrak{n}}| \leqslant \mathcal{N}$ provide $\mathfrak{m}\mathcal{N}$ linearly independent generalised eigenvectors of $\tilde{\mathcal{L}}$ with eigenvalues contained in $A_{\vec{0}}$, these $\mathfrak{m}\mathcal{N}$ eigenvectors must span the centre subspace $\mathbb{E}_c^\mathcal{N}$.

**Lemma 5.** *For generalised eigenvector $\tilde{\mathcal{V}}^{\vec{\mathfrak{n}}}$ constructed from recurrence relations* (28), *derivatives of $\tilde{\mathcal{V}}^{\vec{\mathfrak{n}}}$ with respect to $\vec{\xi}$ satisfy $\partial_{\vec{\xi}}^{\vec{\mathfrak{m}}} \tilde{\mathcal{V}}^{\vec{\mathfrak{n}}} = \tilde{\mathcal{V}}^{\vec{\mathfrak{n}}-\vec{\mathfrak{m}}}$ for $\vec{0} < \vec{\mathfrak{m}} \leqslant \vec{\mathfrak{n}}$.*

Since $\tilde{\mathcal{V}}^{\vec{\mathfrak{n}}} = \frac{\vec{\xi}^{\vec{\mathfrak{n}}}}{\vec{\mathfrak{n}}!}V^{\vec{0}} + \tilde{V}^{\vec{\mathfrak{n}}}$ with $\tilde{V}^{\vec{\mathfrak{n}}}$ of order less than $\vec{\mathfrak{n}}$ in $\vec{\xi}$, then $\partial_{\vec{\xi}}^{\vec{\mathfrak{n}}} \tilde{\mathcal{V}}^{\vec{\mathfrak{n}}} = \tilde{\mathcal{V}}^{\vec{0}}$, in agreement with Lemma 5 when $\vec{\mathfrak{m}} = \vec{\mathfrak{n}}$. However, to prove Lemma 5 we need only prove the $|\vec{\mathfrak{m}}| = 1$ case for general $\vec{\mathfrak{n}}$, as the additive property of derivative powers[3] then ensure the Lemma is true for all $\vec{0} < \vec{\mathfrak{m}} \leqslant \vec{\mathfrak{n}}$.

---
[3] That is, $\partial_{\vec{\xi}}^{\vec{\mathfrak{k}}+\vec{\mathfrak{m}}} \tilde{\mathcal{V}}^{\vec{\mathfrak{n}}} = \partial_{\vec{\xi}}^{\vec{\mathfrak{k}}} \partial_{\vec{\xi}}^{\vec{\mathfrak{m}}} \tilde{\mathcal{V}}^{\vec{\mathfrak{n}}}$ for $\tilde{\mathcal{V}}^{\vec{\mathfrak{n}}}$ a multinomial in $\vec{\xi}$.

We prove Lemma 5 by induction. For $|\vec{n}| = 1$, since $\tilde{\mathcal{V}}^{\vec{n}} = \frac{\vec{\xi}^{\vec{n}}}{\vec{n}!}V^{\vec{0}} + \tilde{V}^{\vec{n}}$ with $\tilde{V}^{\vec{n}}$ independent of $\vec{\xi}$, we know that $\partial_{\vec{\xi}}^{\vec{n}}\tilde{\mathcal{V}}^{\vec{n}} = V^{\vec{0}} = \tilde{\mathcal{V}}^{\vec{0}}$, and thus Lemma 5 is true for all $|\vec{n}| = 1$. Now assume that $\partial_{\vec{\xi}}^{\vec{m}}\tilde{\mathcal{V}}^{\vec{k}} = \tilde{\mathcal{V}}^{\vec{k}-\vec{m}}$ with $|\vec{m}| = 1$ is true for all $\vec{0} < \vec{k} \leqslant \vec{n}$. Then, for $|\vec{m}| = 1$, replace $\vec{n}$ with $\vec{n}+\vec{m}$ in (28b) and take the $\vec{m}$th derivative with respect to $\vec{\xi}$,

$$\begin{aligned}
&\mathfrak{L}_{\vec{0}}(\partial_{\vec{\xi}}^{\vec{m}}\tilde{\mathcal{V}}^{\vec{n}+\vec{m}}) - (\partial_{\vec{\xi}}^{\vec{m}}\tilde{\mathcal{V}}^{\vec{n}+\vec{m}})A_{\vec{0}} \\
&= -\sum_{0<|\vec{k}|,\vec{k}\leqslant\vec{n}+\vec{m}} \mathfrak{L}_{\vec{k}}(\partial_{\vec{\xi}}^{\vec{m}}\tilde{\mathcal{V}}^{\vec{n}+\vec{m}-\vec{k}}) + \sum_{0<|\vec{k}|,\vec{k}\leqslant\vec{n}+\vec{m}} (\partial_{\vec{\xi}}^{\vec{m}}\tilde{\mathcal{V}}^{\vec{n}+\vec{m}-\vec{k}})A_{\vec{k}} \\
&= -\sum_{0<|\vec{k}|,\vec{k}\leqslant\vec{n}} \mathfrak{L}_{\vec{k}}(\partial_{\vec{\xi}}^{\vec{m}}\tilde{\mathcal{V}}^{\vec{n}+\vec{m}-\vec{k}}) + \sum_{0<|\vec{k}|,\vec{k}\leqslant\vec{n}} (\partial_{\vec{\xi}}^{\vec{m}}\tilde{\mathcal{V}}^{\vec{n}+\vec{m}-\vec{k}})A_{\vec{k}} \\
&= -\sum_{0<|\vec{k}|,\vec{k}\leqslant\vec{n}} \mathfrak{L}_{\vec{k}}\tilde{\mathcal{V}}^{\vec{n}-\vec{k}} + \sum_{0<|\vec{k}|,\vec{k}\leqslant\vec{n}} \tilde{\mathcal{V}}^{\vec{n}-\vec{k}}A_{\vec{k}} \\
&= \mathfrak{L}_{\vec{0}}\tilde{\mathcal{V}}^{\vec{n}} - \tilde{\mathcal{V}}^{\vec{n}}A_{\vec{0}},
\end{aligned}$$

where in the third line we recall that the highest $\vec{\xi}$ order of $\tilde{\mathcal{V}}^{\vec{n}+\vec{m}-\vec{k}}$ is $\vec{n}+\vec{m}-\vec{k}$, so to take the $\vec{m}$th derivative we must have $\vec{k} \leqslant \vec{n}$; in the fourth line we apply the assumption that $\partial_{\vec{\xi}}^{\vec{m}}\tilde{\mathcal{V}}^{\vec{k}} = \tilde{\mathcal{V}}^{\vec{k}-\vec{m}}$ for all $\vec{0} < \vec{k} \leqslant \vec{n}$; and the fifth line comes from equation (28b). On comparing the first and last lines we see that $\partial_{\vec{\xi}}^{\vec{m}}\tilde{\mathcal{V}}^{\vec{n}+\vec{m}} - \tilde{\mathcal{V}}^{\vec{n}} \propto V^{\vec{0}}$, but since

$$\begin{aligned}
\partial_{\vec{\xi}}^{\vec{m}}\tilde{\mathcal{V}}^{\vec{n}+\vec{m}} - \tilde{\mathcal{V}}^{\vec{n}} &= \partial_{\vec{\xi}}^{\vec{m}}\left(\frac{\vec{\xi}^{(\vec{n}+\vec{m})}}{(\vec{n}+\vec{m})!}V^{\vec{0}} + \tilde{V}^{\vec{n}+\vec{m}}\right) - \frac{\vec{\xi}^{\vec{n}}}{\vec{n}!}V^{\vec{0}} - \tilde{V}^{\vec{n}} \\
&= \partial_{\vec{\xi}}^{\vec{m}}\tilde{V}^{\vec{n}+\vec{m}} - \tilde{V}^{\vec{n}},
\end{aligned}$$

and we know $\langle Z^{\vec{0}}, \tilde{V}^{\vec{n}}\rangle = 0$ for all $|\vec{n}| > 0$, then $\partial_{\vec{\xi}}^{\vec{m}}\tilde{\mathcal{V}}^{\vec{n}+\vec{m}} = \tilde{\mathcal{V}}^{\vec{n}}$. So, if $\partial_{\vec{\xi}}^{\vec{m}}\tilde{\mathcal{V}}^{\vec{k}} = \tilde{\mathcal{V}}^{\vec{k}-\vec{m}}$ with $|\vec{m}| = 1$ is true for all $\vec{0} < \vec{k} \leqslant \vec{n}$ and $1 \leqslant |\vec{n}| < N$, then $\partial_{\vec{\xi}}^{\vec{m}}\tilde{\mathcal{V}}^{\vec{n}+\vec{m}} = \tilde{\mathcal{V}}^{\vec{n}}$ is also true. Since $\partial_{\vec{\xi}}^{\vec{n}}\tilde{\mathcal{V}}^{\vec{n}} = \tilde{\mathcal{V}}^{\vec{0}}$ is true when $|\vec{n}| = 1$, $\partial_{\vec{\xi}}^{\vec{m}}\tilde{\mathcal{V}}^{\vec{n}+\vec{m}} = \tilde{\mathcal{V}}^{\vec{n}}$ must be true for all $|\vec{n}| > 0$ when $|\vec{m}| = 1$. Finally, because derivative orders are additive, Lemma 5 must be true for all $\vec{0} < \vec{m} \leqslant \vec{n}$.

## 3.2 Slow field and PDE

In this section we complete our primary aim, which is to model the slow dynamics of the original field $u(\vec{x}, y, t)$. To do this, we project $u(\vec{x}, y, t)$ onto the centre subspace $\mathbb{E}_c^0$ and define this projection as the slow field $U(\vec{x}, t) = \langle Z^{\vec{0}}, u(\vec{x}, y, t)\rangle \in \mathbb{C}^m$. The aim of this section is to construct a PDE for $U(\vec{x}, t)$ with an exact error term. For the shallow fluid flow example of Section 1.1, PDEs of different order are (4a) and (5a), but the error term is new.

The slow field $U(\vec{x}, t)$ evaluated at station $\vec{x} = \vec{X}$ is equivalent to $\langle Z^{\vec{0}}, u(\vec{X} + \vec{\xi}, y, t)\rangle$ evaluated at $\vec{\xi} = \vec{0}$, and since $\tilde{u}(\vec{X}, t) = u(\vec{X} + \vec{\xi}, y, t) + \mathcal{O}(|\vec{\xi}|^{N+1})$ (Section 2.2 and equation (23)) the slow field is also equivalent to $\langle Z^{\vec{0}}, \tilde{u}(\vec{X}, t)\rangle$ evaluated at $\vec{\xi} = \vec{0}$. We expand $\tilde{u}(\vec{X}, t)$ in terms of the centre modes $\tilde{\mathcal{V}}^{\vec{n}}$ and the analogous stable modes, and then project this parameterisation onto the the centre subspace $\mathbb{E}_c^0$ to obtain the slow field $U(\vec{X}, t)$.

Since we project $\tilde{u}(\vec{X}, t)$ onto $\mathbb{E}_c^0$, the stable modes may at first seem superfluous in the expansion of $\tilde{u}(\vec{X}, t)$. However, while the stable modes decay exponentially rapidly, resulting in the emergence of the evolution of $u(\vec{x}, y, t)$ on the centre subspace, through the nonlinearity these stable modes are not generally negligible and their influence must be accounted for in $U(\vec{X}, t)$.

We define $\tilde{\mathcal{W}}^{\vec{n}} \in \mathbb{U}^{1 \times m'} \otimes_t \mathbb{G}_N = \mathbb{U}_N^{1 \times m'}$ as the generalised eigenvectors which span the stable subspace $\mathbb{E}_s^N \subset \mathbb{U}_N$ of $\tilde{\mathcal{L}}$. The full set of generalised eigenvectors, $\tilde{\mathcal{V}}^{\vec{n}}$ and $\tilde{\mathcal{W}}^{\vec{n}}$, span $\mathbb{U}_N$ of $\tilde{\mathcal{L}}$, fully parameterising the field $\tilde{u}$. Many of the properties of the $\tilde{\mathcal{W}}^{\vec{n}}$ are analogous to those of the $\tilde{\mathcal{V}}^{\vec{n}}$ and can be established by proofs similar to those presented in Section 3.1. Therefore, here we only briefly comment on those properties of $\tilde{\mathcal{W}}^{\vec{n}}$ which are required for the analysis of generating multinomial $\tilde{u}$, and ultimately the slow field $U(\vec{x}, t)$.

For the lowest order case $\vec{n} = \vec{0}$, the centre subspace eigenvector $\tilde{\mathcal{V}}^{\vec{0}} = V^{\vec{0}}$ satisfies $\mathfrak{L}_{\vec{0}} \tilde{\mathcal{V}}^{\vec{0}} = \tilde{\mathcal{V}}^{\vec{0}} A_{\vec{0}}$ (Definition 4) where matrix $A_{\vec{0}} \in \mathbb{C}^{m \times m}$ has centre eigenvalues $\lambda_1, \ldots, \lambda_m$ (Assumption 3), and similarly there is a stable subspace eigenvector $\tilde{\mathcal{W}}^{\vec{0}} = \begin{bmatrix} v_{m+1}^{\vec{0}} & v_{m+2}^{\vec{0}} & \cdots \end{bmatrix} \in \mathbb{U}^{1 \times m'}$ which satisfies $\mathfrak{L}_{\vec{0}} \tilde{\mathcal{W}}^{\vec{0}} = \tilde{\mathcal{W}}^{\vec{0}} B_{\vec{0}}$ for some matrix $B_{\vec{0}} \in \mathbb{C}^{m' \times m'}$ with stable eigenvalues $\lambda_{m+1}, \lambda_{m+2}, \ldots, \lambda_{m'}$ (Recall from Assumption 3 that $m'$ may be infinite). Together, $\tilde{\mathcal{V}}^{\vec{0}}$ and $\tilde{\mathcal{W}}^{\vec{0}}$ span the Hilbert space $\mathbb{U}$ of $\mathfrak{L}_{\vec{0}}$.

For convenience, define matrix $\tilde{\mathcal{V}} = [\tilde{\mathcal{V}}^{\vec{n}}]$ where the columns of $\tilde{\mathcal{V}}$ are the centre subspace eigenvectors $\tilde{\mathcal{V}}^{\vec{n}}$. The ordering of the $\tilde{\mathcal{V}}^{\vec{n}}$ in $\tilde{\mathcal{V}}$ is according to the magnitude $|\vec{n}|$, so that the first column is $\tilde{\mathcal{V}}^{\vec{0}}$. From (29), $\tilde{\mathcal{L}} \tilde{\mathcal{V}}^{\vec{n}} = \sum_{\vec{0} \leqslant \vec{k} \leqslant \vec{n}} \tilde{\mathcal{V}}^{\vec{n}-\vec{k}} A_{\vec{k}}$ and so we define block upper triangular matrix $\mathcal{A}$ such that $\tilde{\mathcal{L}} \tilde{\mathcal{V}} = \tilde{\mathcal{V}} \mathcal{A}$. The upper block triangular matrix $\mathcal{A}$ consists of $\mathcal{N} \times \mathcal{N}$ blocks with $0_m$ below the diagonal, $A_{\vec{0}} \in \mathbb{C}^{m \times m}$ along the main diagonal, and the $(\vec{k}, \vec{n})$ block above the diagonal (that is, for $\vec{n} > \vec{k}$) is $A_{\vec{n}-\vec{k}} \in \mathbb{C}^{m \times m}$. The centre eigenvalues of $\tilde{\mathcal{L}}$ are the eigenvalues of the blocks along the diagonal of $\mathcal{A}$, namely, the $m$ eigenvalues of $A_{\vec{0}}$, $\lambda_1, \ldots, \lambda_m$, repeated $\mathcal{N}$ times.

Similarly to matrix $\mathcal{A}$ which satisfies $\tilde{\mathcal{L}} \tilde{\mathcal{V}} = \tilde{\mathcal{V}} \mathcal{A}$, we define matrix $\mathcal{B}$ such that $\tilde{\mathcal{L}} \tilde{\mathcal{W}} = \tilde{\mathcal{W}} \mathcal{B}$. Analogous to $\mathcal{A}$, $\mathcal{B}$ is upper block triangular with $\mathcal{N} \times \mathcal{N}$ blocks of size $\mathbb{C}^{m' \times m'}$, with block $B_{\vec{0}}$ along the main diagonal, and $B_{\vec{n}-\vec{k}}$ the $(\vec{k}, \vec{n})$ block above the diagonal. Thus the stable eigenvalues of $\tilde{\mathcal{L}}$ are the eigenvalues of $\mathcal{B}$ which must be $\lambda_{m+1}, \lambda_{m+2}, \ldots$ repeated $\mathcal{N}$ times. Recall that these stable eigenvalues all have real part $\leqslant -\beta < -N\alpha$, whereas the magnitude of the real part of the centre eigenvalues are $\leqslant \alpha$ (Assumption 3).

We capture the full centre-stable dynamics of the linear generating PDE (27) on $\mathbb{U}_N$ with

$$\tilde{u}(\vec{X}, t) = \sum_{|\vec{n}|=0}^{N} (\tilde{\mathcal{V}}^{\vec{n}} U_{\vec{n}} + \tilde{\mathcal{W}}^{\vec{n}} S_{\vec{n}}) = \tilde{\mathcal{V}} \mathcal{U} + \tilde{\mathcal{W}} \mathcal{S}, \quad (30)$$

for parameters $\mathcal{U} = [U_{\vec{n}}]$ and $\mathcal{S} = [S_{\vec{n}}]$ with $U_{\vec{n}} \in \mathbb{C}^m$ and $S_{\vec{n}} \in \mathbb{C}^{m'}$. As $\tilde{r}[u], \tilde{f}[\tilde{u}] \in \mathbb{U}_N$ and the generalised eigenvectors, $\tilde{\mathcal{V}}^{\vec{n}}$ and $\tilde{\mathcal{W}}^{\vec{n}}$, span $\mathbb{U}_N$ the forcing and nonlinear terms are uniquely parameterised in terms of these generalised eigenvectors,

$$\tilde{r}[u] = \tilde{\mathcal{V}} r_c(t) + \tilde{\mathcal{W}} r_s(t), \quad \tilde{f}[\tilde{u}] = \tilde{\mathcal{V}} f_c(t) + \tilde{\mathcal{W}} f_s(t), \quad (31)$$

where $r_c = [r_c^{\vec{n}}]$ and $f_c = [f_c^{\vec{n}}]$, and similarly for $r_s$ and $f_s$, with $r_c^{\vec{n}}, f_c^{\vec{n}} \in \mathbb{C}^m$ and $r_s^{\vec{n}}, f_s^{\vec{n}} \in \mathbb{C}^{m'}$. We substitute the expansion of $\tilde{u}$ (30) into the nonlinear generating PDE (18) and separate the forcing and nonlinear terms into centre and slow components (31). From Section 3.1 we know that the centre subspace eigenvectors $\tilde{\mathcal{V}}^{\vec{n}}$ are linearly independent, and similarly the stable subspace eigenvectors $\tilde{\mathcal{W}}^{\vec{n}}$ are linearly independent. Therefore, we separate the centre and slow components of

the PDE to obtain

$$\frac{\partial \mathcal{U}}{\partial t} = \mathcal{A}\mathcal{U} + f_c(t) + r_c(t), \tag{32a}$$

$$\frac{\partial \mathcal{S}}{\partial t} = \mathcal{B}\mathcal{S} + f_s(t) + r_s(t). \tag{32b}$$

A general solution of the PDE for the stable parameter $\mathcal{S}$ (32b) is

$$\mathcal{S}(t) = e^{\mathcal{B}t}\mathcal{S}(0) + \int_0^t e^{\mathcal{B}(t-\tau)}[f_s(\tau) + r_s(\tau)]d\tau = e^{\mathcal{B}t}\mathcal{S}(0) + e^{\mathcal{B}t} \star [f_s(t) + r_s(t)], \tag{33}$$

with convolution $h(t) \star g(t) := \int_0^t h(t-\tau)g(\tau)\,d\tau$. As all eigenvalues of $\mathcal{B}$ have real part $\leqslant -\beta < -N\alpha$, for some decay rate $\gamma \in (\alpha, \beta)$,

$$\mathcal{S}(t) = e^{\mathcal{B}t} \star [f_s(t) + r_s(t)] + \mathcal{O}\left(e^{-\gamma t}\right). \tag{34}$$

This solution for $\mathcal{S}$ shows that, after a sufficiently long time, the forcing and nonlinear terms dominate $\mathcal{S}$ through the convolution, thus showing how the forcing and nonlinear terms couple the centre and stable solutions through $f_{s/c}$ and $r_{s/c}$ and why the influence of the stable modes are not negligible.

We now construct a PDE for the slow field $U(\vec{X}, t)$ by considering both its the temporal and spatial derivative in terms of the centre-stable dynamics. Since $U(\vec{X}, t) = \langle Z^{\vec{0}}, \tilde{u}(\vec{X}, t)\rangle_{\vec{\xi}=\vec{0}}$ at station $\vec{x} = \vec{X}$,

$$\begin{aligned}
\frac{\partial U(\vec{X}, t)}{\partial t} &= \left\langle Z^{\vec{0}}, \frac{\partial \tilde{u}(\vec{X}, t)}{\partial t}\right\rangle_{\vec{\xi}=\vec{0}} \\
&= \left\langle Z^{\vec{0}}, \tilde{\mathcal{V}}\frac{\partial \mathcal{U}}{\partial t}\right\rangle_{\vec{\xi}=\vec{0}} + \left\langle Z^{\vec{0}}, \tilde{\mathcal{W}}\frac{\partial \mathcal{S}}{\partial t}\right\rangle_{\vec{\xi}=\vec{0}} \\
&= \langle Z^{\vec{0}}, \tilde{\mathcal{V}}\mathcal{A}\mathcal{U}\rangle_{\vec{\xi}=\vec{0}} + \langle Z^{\vec{0}}, \tilde{\mathcal{V}}f_c(t)\rangle_{\vec{\xi}=\vec{0}} + \langle Z^{\vec{0}}, \tilde{\mathcal{V}}r_c(t)\rangle_{\vec{\xi}=\vec{0}} \\
&\quad + \langle Z^{\vec{0}}, \tilde{\mathcal{W}}\mathcal{B}\mathcal{S}\rangle_{\vec{\xi}=\vec{0}} + \langle Z^{\vec{0}}, \tilde{\mathcal{W}}f_s(t)\rangle_{\vec{\xi}=\vec{0}} + \langle Z^{\vec{0}}, \tilde{\mathcal{W}}r_s(t)\rangle_{\vec{\xi}=\vec{0}} \quad \text{from (32)} \\
&= \langle Z^{\vec{0}}, \tilde{\mathcal{V}}\rangle_{\vec{\xi}=\vec{0}}\mathcal{A}\mathcal{U} + \langle Z^{\vec{0}}, \tilde{\mathcal{V}}\rangle_{\vec{\xi}=\vec{0}}f_c(t) + \langle Z^{\vec{0}}, \tilde{\mathcal{V}}\rangle_{\vec{\xi}=\vec{0}}r_c(t) \\
&\quad + \langle Z^{\vec{0}}, \tilde{\mathcal{W}}\rangle_{\vec{\xi}=\vec{0}}f_s(t) + \langle Z^{\vec{0}}, \tilde{\mathcal{W}}\rangle_{\vec{\xi}=\vec{0}}r_s(t) \\
&\quad + \langle Z^{\vec{0}}, \tilde{\mathcal{W}}\rangle_{\vec{\xi}=\vec{0}}\mathcal{B}e^{\mathcal{B}t} \star [f_s(t) + r_s(t)] + \mathcal{O}\left(e^{-\gamma t}\right) \quad \text{from (34)}.
\end{aligned} \tag{35}$$

From (28c), the inner product $\langle Z^{\vec{0}}, \tilde{\mathcal{V}}\rangle_{\vec{\xi}=0} = \begin{bmatrix} I_m & 0_m & \cdots & 0_m \end{bmatrix}$. Also, because of the upper block triangular structure of $\mathcal{A}$, the $\vec{k}$ element of $\mathcal{A}\mathcal{U}$ is

$$[\mathcal{A}\mathcal{U}]_{\vec{k}} = \sum_{\vec{n} \geqslant \vec{k}, |\vec{n}| \leqslant N} \mathcal{A}_{\vec{n}-\vec{k}} U_{\vec{n}} = \sum_{|\vec{n}|=0}^{N-|\vec{k}|} \mathcal{A}_{\vec{n}} U_{\vec{n}+\vec{k}}.$$

Therefore, in the first term on the right hand side of (35) we only retain the $\vec{k} = \vec{0}$ element of $\mathcal{A}\mathcal{U}$, and in the second and third terms we only retain $f_c^{\vec{0}}$ and $r_c^{\vec{0}}$, respectively. So now,

$$\begin{aligned}
\frac{\partial U(\vec{X}, t)}{\partial t} &= \sum_{|\vec{n}|=0}^{N} \mathcal{A}_{\vec{n}} U_{\vec{n}} + f_c^{\vec{0}}(t) + r_c^{\vec{0}}(t) \\
&\quad + \langle Z^{\vec{0}}, \tilde{\mathcal{W}}\rangle_{\vec{\xi}=\vec{0}}f_s(t) + \langle Z^{\vec{0}}, \tilde{\mathcal{W}}\rangle_{\vec{\xi}=\vec{0}}r_s(t)
\end{aligned}$$

$$+ \langle \vec{Z^0}, \tilde{\mathcal{W}} \rangle_{\vec{\xi}=\vec{0}} \mathcal{B} e^{\mathcal{B}t} \star [f_s(t) + r_s(t)] + \mathcal{O}(e^{-\gamma t}). \quad (36)$$

Now consider the order $\vec{n}$ spatial derivative of the slow field,

$$\begin{aligned} [\partial_{\vec{x}}^{\vec{n}} U(\vec{x}, t)]_{\vec{x}=\vec{X}} &= \langle \vec{Z^0}, u(\vec{x}, y, t) \rangle_{\vec{x}=\vec{X}} \\ &= \langle \vec{Z^0}, u^{(\vec{n})} \rangle \\ &= \langle \vec{Z^0}, \partial_{\vec{\xi}}^{\vec{n}} \tilde{u}(\vec{X}, t) \rangle_{\vec{\xi}=\vec{0}} \quad \text{from (21)} \\ &= \langle \vec{Z^0}, \partial_{\vec{\xi}}^{\vec{n}} \tilde{\mathcal{V}} \rangle_{\vec{\xi}=\vec{0}} \mathcal{U} + \langle \vec{Z^0}, \partial_{\vec{\xi}}^{\vec{n}} \tilde{\mathcal{W}} \rangle_{\vec{\xi}=\vec{0}} \mathcal{S}. \end{aligned} \quad (37)$$

From Lemma 5, $\langle \vec{Z^0}, \partial_{\vec{\xi}}^{\vec{n}} \tilde{\mathcal{V}}^{\vec{k}} \rangle_{\vec{\xi}=\vec{0}} = \langle \vec{Z^0}, \tilde{\mathcal{V}}^{\vec{k}-\vec{n}} \rangle_{\vec{\xi}=\vec{0}}$ which, from (28c), equals the identity $I_m$ if $\vec{k} = \vec{n}$, but zero otherwise. Thus, in the first inner product of (37), only the $U_{\vec{n}}$ element of $\mathcal{U}$ remains. Then, in the second inner product of (37), substitute the solution of the stable parameter $\mathcal{S}$ (34). The spatial derivative is now

$$[\partial_{\vec{x}}^{\vec{n}} U(\vec{x}, t)]_{\vec{x}=\vec{X}} = U_{\vec{n}} + \langle \vec{Z^0}, \partial_{\vec{\xi}}^{\vec{n}} \tilde{\mathcal{W}} \rangle_{\vec{\xi}=\vec{0}} e^{\mathcal{B}t} \star [f_s(t) + r_s(t)] + \mathcal{O}(e^{-\gamma t}). \quad (38)$$

Combining (36) and (38),

$$\begin{aligned} \frac{\partial U(\vec{X}, t)}{\partial t} &= \sum_{|\vec{n}|=0}^{N} A_{\vec{n}} \partial_{\vec{x}}^{\vec{n}} U(\vec{X}, t) \\ &\quad + f_c^{\vec{0}}(t) + r_c^{\vec{0}}(t) + \langle \vec{Z^0}, \tilde{\mathcal{W}} \rangle_{\vec{\xi}=\vec{0}} f_s(t) + \langle \vec{Z^0}, \tilde{\mathcal{W}} \rangle_{\vec{\xi}=\vec{0}} r_s(t) \\ &\quad + \left[ \langle \vec{Z^0}, \mathcal{W} \rangle \mathcal{B} - \sum_{|\vec{n}|=0}^{N} A_{\vec{n}} \langle \vec{Z^0}, \partial_{\vec{\xi}}^{\vec{n}} \tilde{\mathcal{W}} \rangle \right]_{\vec{\xi}=\vec{0}} e^{\mathcal{B}t} \star [f_s(t) + r_s(t)] \\ &\quad + \mathcal{O}(e^{-\gamma t}). \end{aligned} \quad (39)$$

Whereas this equation symbolically resembles a PDE, it is strictly a differential-integral equation which couples the dynamics at each station $\vec{X}$ via the 'uncertain' gradient terms and the stable parameter $\mathcal{S}(t)$, which is dependent on the history convolution integrals (33). To obtain a slow PDE without this coupling to different stations, such as PDE (4a) used in the shallow fluid flow example, we retain all terms which do not couple to different stations (i.e., no dependence on derivatives $\partial_{\vec{x}} u^{(\vec{n})}$ with $|\vec{n}| = N$ and no dependence on $\mathcal{S}(t)$) and regulate all other terms to a *remainder*. The last line of (39) contains a convolution, so is part of the remainder, and the two forcing terms $r_c^{\vec{0}}(t)$ and $\langle \vec{Z^0}, \tilde{\mathcal{W}} \rangle_{\vec{\xi}=\vec{0}} r_s(t)$ are dependent on uncertain gradients, so are also in the remainder. In contrast, the nonlinear terms $f_c^{\vec{0}}(t)$ and $f_s(t)$ contain parts which we want to retain in the slow PDE, as well as terms which should be in the remainder.

For specific cases, removing the remainder components from the nonlinear terms in (39) is achieved using (26), as shown in Appendix A.2. Here, for the general case, we show that the nonlinear terms which are retained in the slow PDE must take a particular form. First, separate the nonlinear terms in the second line of (39) into two parts,

$$f_c^{\vec{0}}(t) = f_c^{\vec{0}}(\vec{X}, t) + f_{c,r}^{\vec{0}}(t), \quad \langle \vec{Z^0}, \tilde{\mathcal{W}} \rangle_{\vec{\xi}=\vec{0}} f_s(t) = f_s(\vec{X}, t) + f_{s,r}(t), \quad (40)$$

where $f_c^{\vec{0}}(\vec{X}, t) \in \mathbb{C}^m$ and $f_s^{\vec{0}}(\vec{X}, t) \in \mathbb{C}^m$ contain no uncertain terms and no dependence on $\mathcal{S}$ (so are retained in the slow PDE), and where $f_{c,r}^{\vec{0}}(t)$ and $f_{s,r}(t)$ contain all uncertain terms and $\mathcal{S}$ dependent terms. The $f^{\vec{0}}_{c,r}(t)$ and $f_{s,r}(t)$, as well as the convolutions and forcing terms in (39), are not in retained in the slow PDE. As the nonlinear function $f[u]$ in the original PDE (6) is a sum of nonlinear terms $f^j[u]$ (Assumption 1), the nonlinear $f^{\vec{0}}_c(\vec{X},t)$ and $f_s(\vec{X},t)$ are also a sum of nonlinear terms indexed by integer $j$,

$$f^{\vec{0}}_c(\vec{X},t) = \sum_j f^{j\vec{0}}_c(\vec{X},t), \quad f^{\vec{0}}_s(\vec{X},t) = \sum_j f^{j\vec{0}}_s(\vec{X},t), \qquad (41)$$

where $f^{j\vec{0}}_c(\vec{X},t) \in \mathbb{C}^m$, and $f^{j\vec{0}}_s(\vec{X},t) \in \mathbb{C}^m$. As $f^j[u]$ is of order $P_j$ in $u$ and its derivatives (Assumption 1), $f^{j\vec{0}}_c$ and $f^{j\vec{0}}_s$ must be of order $P_j$ in $U_{\vec{n}} \in \mathbb{C}^m$ for all $|\vec{n}| \leq N$ [4]. So, in general, each $k = 1, \ldots, m$ element of $f^j_c + f^j_s$ must have the form

$$\sum_{\substack{|\vec{\ell}_1|,\ldots,|\vec{\ell}_{P_j}|=0 \\ |\vec{\ell}_1| \geq |\vec{\ell}_2| \geq \cdots \geq |\vec{\ell}_{P_j}|}}^{N} \vec{a}^{jT}_{k\vec{\ell}_1\vec{\ell}_2\ldots\vec{\ell}_{P_j}} U_{\vec{\ell}_1} \otimes U_{\vec{\ell}_2} \otimes \cdots \otimes U_{\vec{\ell}_{P_j}}, \qquad (42)$$

for some constant vector $\vec{a}^j_{k\vec{\ell}_1\vec{\ell}_2\ldots\vec{\ell}_{P_j}} \in \mathbb{C}^{m^{P_j}}$ and where $\otimes$ represents the usual Kronecker product, for which $U_{\vec{\ell}_p} \otimes U_{\vec{\ell}_q} \in \mathbb{C}^{m^2}$ and $U_{\vec{\ell}_1} \otimes U_{\vec{\ell}_2} \otimes \cdots \otimes U_{\vec{\ell}_{P_j}} \in \mathbb{C}^{m^{P_j}}$. [5] On replacing all $U_{\vec{\ell}}$ with spatial derivatives of $\partial^{\vec{\ell}}_{\vec{x}} U$, as shown in (38), the $k$th coordinate of the $m$-dimensional nonlinear term retained in the slow PDE must have the form

$$[f^{\vec{0}}_c(\vec{X},t) + f_s(\vec{X},t)]_k$$
$$= \sum_j \sum_{\substack{|\vec{\ell}_1|,\ldots,|\vec{\ell}_{P_j}|=0 \\ |\vec{\ell}_1| \geq |\vec{\ell}_2| \geq \cdots \geq |\vec{\ell}_{P_j}|}}^{N} \vec{a}^{jT}_{k\vec{\ell}_1\vec{\ell}_2\ldots\vec{\ell}_{P_j}} (\partial^{\vec{\ell}_1}_{\vec{x}} U) \otimes (\partial^{\vec{\ell}_2}_{\vec{x}} U) \otimes \cdots \otimes (\partial^{\vec{\ell}_{P_j}}_{\vec{x}} U). \qquad (43)$$

Now, on replacing arbitrary station $\vec{X}$ with $\vec{x} \in \mathbb{X}$, the slow PDE determined from the differential-integral equation (39) is

$$\frac{\partial U(\vec{x},t)}{\partial t} = \sum_{|\vec{n}|=0}^{N} A_{\vec{n}} \partial^{\vec{n}}_{\vec{x}} U(\vec{x},t)$$
$$+ \sum_j \sum_{\substack{|\vec{\ell}_1|,\ldots,|\vec{\ell}_{P_j}|=0 \\ |\vec{\ell}_1| \geq |\vec{\ell}_2| \geq \cdots \geq |\vec{\ell}_{P_j}|}}^{N} \vec{a}^{jT}_{\vec{\ell}_1\vec{\ell}_2\ldots\vec{\ell}_{P_j}} (\partial^{\vec{\ell}_1}_{\vec{x}} U) \otimes (\partial^{\vec{\ell}_2}_{\vec{x}} U) \otimes \cdots \otimes (\partial^{\vec{\ell}_{P_j}}_{\vec{x}} U) + \rho, \qquad (44)$$

with $\vec{a}^{jT}_{\vec{\ell}_1\vec{\ell}_2\ldots\vec{\ell}_{P_j}} (\partial^{\vec{\ell}_1}_{\vec{x}} U) \otimes \cdots \otimes (\partial^{\vec{\ell}_{P_j}}_{\vec{x}} U)$ the $m$-dimensional vector with elements $k = 1, 2, \ldots, m$ defined by (43), and with remainder

$$\rho = f^{\vec{0}}_{c,r}(t) + r^{\vec{0}}_c(t) + f_{s,r}(t) + \langle Z^{\vec{0}}, \tilde{\mathcal{W}} \rangle_{\vec{\xi}=\vec{0}} r_s(t)$$

---

[4] The nonlinear $f^{\vec{0}}_c(t)$ and $\langle Z^{\vec{0}}, \tilde{\mathcal{W}} \rangle_{\vec{\xi}=\vec{0}} f_s(t)$ are sums of nonlinear terms of order $P_j$ in both $U_{\vec{n}}$ and $S_{\vec{n}}$ for all $|\vec{n}| \leq N$, but as all $S_{\vec{n}}$ dependence is contained in $f^{\vec{0}}_{c,r}(t)$ and $f_{s,r}(t)$ there is only $U_{\vec{n}}$ dependence in $f^{\vec{0}}_c(\vec{X},t)$ and $f_s(\vec{X},t)$.

[5] For $m = 1$, as in the fluid flow example of Section 1.1, the constant vector in (42) reduces to a scalar and the Kronecker products reduce to a multiplication of scalars. For $P_j = 2$ and any value of $m$, each term in (42) is equivalent to $U^T_{\vec{\ell}_1} C^j_{k\vec{\ell}_1\vec{\ell}_2} U_{\vec{\ell}_2}$ for some constant matrix $C^j_{k\vec{\ell}_1\vec{\ell}_2} \in \mathbb{C}^{m \times m}$.

$$+ \left[ \langle Z^{\vec{0}}, \tilde{\mathcal{W}} \rangle \mathcal{B} - \sum_{|\vec{n}|=0}^{N} A_{\vec{n}} \langle Z^{\vec{0}}, \partial_{\vec{\xi}}^{\vec{n}} \tilde{\mathcal{W}} \rangle \right]_{\vec{\xi}=\vec{0}} e^{\mathcal{B}t} \star [f_s(t) + r_s(t)] + \mathcal{O}\left(e^{-\gamma t}\right).$$
(45)

Analogous slow PDEs were derived by Roberts (2015) (equation (22)) and Roberts and Bunder (2017) (equation (51)), but without the nonlinear terms.

Simplifications of the remainder $\rho$ (45) are possible when the order $N$ is chosen to be higher than the order of the spatial derivatives in the original PDE (6). The original PDE contains linear operators $\mathfrak{L}_{\vec{k}} \partial_{\vec{x}}^{\vec{k}}$ for $\vec{k}$ satisfying $0 \leqslant |\vec{k}| < \infty$, but in practice there will be an upper limit on $|\vec{k}|$, say $k_{max}$ (often $k_{max} = 1, 2$—the example of Section 1.1 has $k_{max} = 1$). Assume that $N > k_{max}$ and consider the uncertain linear terms in $\rho$:

$$r_c^{\vec{0}}(t) + \langle Z^{\vec{0}}, \tilde{\mathcal{W}} \rangle_{\vec{\xi}=\vec{0}} r_s(t) = \langle Z^{\vec{0}}, \tilde{\mathcal{V}} r_c(t) + \tilde{\mathcal{W}} r_s(t) \rangle_{\vec{\xi}=\vec{0}} = \langle Z^{\vec{0}}, \tilde{r}[u] \rangle_{\vec{\xi}=\vec{0}}.$$

Since $Z^{\vec{0}}$ is independent of $\vec{\xi}$, we need only consider $\vec{\xi} = \vec{0}$ in $\tilde{r}[u]$ (22). When $\vec{\xi} = \vec{0}$ the right hand side of (22) requires $|\vec{n}| = 0$, $|\vec{\ell}| = N$ and $\vec{\ell} \lneq \vec{k}$, but since $N > k_{max} \geqslant |\vec{k}|$ we can never satisfy $|\vec{\ell}| = N$ and $\vec{\ell} \lneq \vec{k}$. So, when $N > k_{max}$ we have $\tilde{r}[u]_{\vec{\xi}=0} = 0$ and $r_c^{\vec{0}}(t) + \langle Z^{\vec{0}}, \tilde{\mathcal{W}} \rangle_{\vec{\xi}=\vec{0}} r_s(t) = 0$. Similarly, consider the projection of the nonlinear term

$$f_c^{\vec{0}}(t) + \langle Z^{\vec{0}}, \tilde{\mathcal{W}} \rangle_{\vec{\xi}=\vec{0}} f_s(t) = \langle Z^{\vec{0}}, \tilde{\mathcal{V}} f_c(t) + \tilde{\mathcal{W}} f_s(t) \rangle_{\vec{\xi}=\vec{0}} = \langle Z^{\vec{0}}, \tilde{f}[\tilde{u}] \rangle_{\vec{\xi}=\vec{0}}.$$

and then expand $\tilde{f}[\tilde{u}]$ using (25). If $N$ is chosen to be larger than any spatial derivative in the nonlinear term, that is $N > p_i^j$ for all $i = 1, 2, \ldots, P_j$ and for all $j$ the number of nonlinear terms, then

$$\tilde{f}[\tilde{u}]_{\vec{\xi}=\vec{0}} = \sum_j c_j(y) \prod_{i=1}^{P_j} \left[ \partial_{\vec{\xi}}^{\vec{p}_i^j} \tilde{u} \right]_{\vec{\xi}=\vec{0}}$$

contains no uncertain terms. So, when separating the nonlinear terms according to (31) and (40) $f_{c,r}^{\vec{0}}(t)$ and $f_{s,r}(t)$ contain all $\mathcal{S}$ dependence and any convolution terms, but no uncertain terms.

We have shown that $N > \max(p_i^j, k_{max})$ removes the uncertain terms from $f_{c,r}^{\vec{0}}(t)$ and $f_{s,r}(t)$ and sets $r_c^{\vec{0}}(t) + \langle Z^{\vec{0}}, \tilde{\mathcal{W}} \rangle_{\vec{\xi}=\vec{0}} r_s(t) = 0$, but this does not remove all uncertain terms from the remainder $\rho$ (45). Uncertain terms are still present in the remainder because of $f_s(t)$ and $r_s(t)$ which appear in the convolution in the second line of (45).

Section 1.1 presents the example of a shallow fluid flow on a rotating substrate and, with computer algebra code provided in Appendix A, constructs slow PDE of the form given in (44), as shown in equations (4a) and (5a), for $N = 3$ and $N = 4$, respectively. Appendix A.3 calculates parts of the remainder $\rho$, such as $\langle Z^{\vec{0}}, \tilde{\mathcal{W}} \rangle \mathcal{B}$ and $A_{\vec{n}} \langle Z^{\vec{0}}, \partial_{\vec{\xi}}^{\vec{n}} \tilde{\mathcal{W}} \rangle$, and shows that, since the order is sufficiently large ($N > \max(p_i^j, k_{max}) = 1$) we have $r_c^{\vec{0}}(t) + \langle Z^{\vec{0}}, \tilde{\mathcal{W}} \rangle_{\vec{\xi}=\vec{0}} r_s(t) = 0$. Whereas the appendix is written to support the example presented in Section 1.1, only Appendix A.1 is specific to this example, with the code in Appendices A.2 and A.3 written in a general format so as to be readily adaptable to a large number of systems.

# 4 Conclusion

This article further develops a general theory to support practical approximations of slow variations in space. This methodology was initially developed by Roberts (2015) for one dimensional space, and then extended by Roberts and Bunder (2017) to linear systems in multi-dimensional space. We here provide theoretical support and a practical example for the case of a nonlinear system of PDEs in a spatial domain that is large in multiple dimensions. The significant advantages of the theoretical methodology are:

- the approach is readily applicable to a wide range of systems, as illustrated by the general theory provided in Sections 2 and 3;

- higher order PDE are obtained in a straightforward manner by increasing the order $N$ of the Taylor expansion; and

- the resulting slow PDE has a well-defined error, with a derived an algebraic form which can be bounded in applications.

In the general theory, we make some assumptions about the structure of the nonlinear microscale system and its dynamics. Assumption 1 requires that the nonlinearity in the original microscale PDE should be a sum of products of the unknown field and its derivatives, and Assumption 3 requires centre-stable dynamics. The key requirement for the presented methodology is the persistence of the centre manifold of the linear system (described in Section 3.1) when perturbed by nonlinearities and time-dependent forcing, thus justifying the importance of the eigenspace of the linearised system to the full nonlinear microscale system (Section 3.2). As other invariant manifolds are often similarly persistent, we expect that the methodology is not restricted to the centre-stable dynamics required by Assumption 3. Indeed, we show that other invariant manifolds are possible with the fluid flow example in Section 1.1, which has slow-stable dynamics. Furthermore, this fluid flow example has a more complex nonlinear structure than that required by Assumption 1. The Reduce Algebra code presented in Appendix A is designed for this fluid flow example, but is written so as to be adaptable to other systems, including those with different nonlinear structures.

Future research will aim to further generalise the methodology. Of particular interest is stochastic dynamics (Arnold and Imkeller 1998; Roberts 2008), deriving boundary conditions for the slowly varying model from microscale boundary conditions (Segel 1969; Roberts 1992; Mielke 1992), and non-local operators (i.e., beyond the local operators $\mathfrak{L}_{\vec{k}}\partial_{\vec{x}}^{\vec{k}}$ in (6)) (Calcagni, Montobbio, and Nardelli 2008).

**Acknowledgement**   The Australian Research Council Discovery Project grants DP150102385 and DP180100050 helped support this research. We thank Arthur Norman and colleagues who maintain the Reduce software.

# A Computer algebra determines the emergent macroscale model

The Reduce Algebra code presented here takes the microscale fluid flow equations (1) for field $u(\vec{x}, y, t)$ and derives the emergent slow macroscale dynamics (4a) or (5a) for the slow field $U(\vec{x}, t)$ on $\vec{x} \in \mathbb{X}$, where domain $\mathbb{X} \subset \mathbb{R}^M$ with $M = 2$ is specified as 'large'. The microscale system (1) decomposes into a stable subspace and an $m = 1$ dimensional slow subspace, with the slow field evolving on this slow subspace $U \in \mathbb{R}^m$.

The Reduce code is applicable to other systems with original microscale field $u \in \mathbb{U} \subset \mathbb{R}^d$ on an $M = 2$ dimensional large domain, and described by PDEs with an order $P = 2$ nonlinear term, provided the dynamics of these systems decompose into a $d-1$ stable subspace and an $m = 1$ dimensional centre or slow subspace, thus producing a slow field $U(\vec{x}, t) \in \mathbb{R}$. But this code is readily adaptable to systems with different specifications. Other parameters in the code are easily changed; for example, the dimension $d$ of the field $u$ and the order of the Taylor expansion $N$ are variables in the computer algebra.

We firstly set some printing options and ensure that Reduce will provide us with complex number solutions. We then specify the dimensions of the $u$ field and the desired order of the Taylor expansion $N$.

```
1  on div; off allfac; on revpri; on complex;
2  d:=3;   % dimension of u field
3  nn:=3; % order of Taylor expansion
```

Appendix A.1 defines the microscale PDE of field $u(\vec{x}, y, t)$, introduced in Section 1.1, which describe fluid flow on a rotating substrate. Appendix A.2 is generic code for any system with an $m = 1$ dimensional slow subspace, a 'large' domain of dimension $M = 2$, and nonlinearity of order $P = 2$. This generic code constructs both slow eigenvectors and eigenvalues, $\tilde{\mathcal{V}}^{\vec{n}}$ and $A_{\vec{n}}$, and stable eigenvectors and eigenvalues, $\tilde{\mathcal{W}}^{\vec{n}}$ and $B_{\vec{n}}$, of the matrix operator $\tilde{\mathcal{L}}$ (19) for $|\vec{n}| \leqslant N$. Then the code produces the slow macroscale PDE for the slow field $U(\vec{x}, t) \in \mathbb{R}^m$ with $\vec{x}$ defined on the large domain $\mathbb{X} \subset \mathbb{R}^M$.

## A.1 Thin film flow

The provided code allows for a more general PDE than that given in (1), by including the kinematic viscosity $\nu$, and Weber number We to describe surface tension,

$$\frac{\partial h}{\partial t} = -\nabla \cdot (h\vec{v}), \tag{46a}$$

$$\frac{\partial \vec{v}}{\partial t} = \begin{bmatrix} -b & f \\ -f & -b \end{bmatrix} \vec{v} - (\vec{v} \cdot \nabla)\vec{v} - g\nabla h + \nu \nabla^2 \vec{v} + \text{We}\, \nabla^3 \vec{v}, \tag{46b}$$

although here, for simplicity, we set $\nu, \text{We} = 0$.

The code below defines matrices $\mathfrak{L}_{\vec{n}}$ and $\mathfrak{M}_{\vec{n},j}$ for $|\vec{n}| \leqslant N$ for PDEs (1) and (46).

```
4  nu:=0$
5  we:=0$ % simple example
6
7  % linear L matrices
8  l100:=mat((0, 0, 0),(0, -b, f), (0, -f, -b));
9  l110:=mat((0, -1, 0),(-g, 0, 0), (0, 0, 0));
10 l101:=mat((0, 0, -1),(0, 0, 0), (-g, 0, 0));
```

```
11 l120:=nu*mat((0, 0, 0),(0, 1, 0), (0, 0, 1));
12 l102:=nu*mat((0, 0, 0),(0, 1, 0), (0, 0, 1));
13 l130:=we*mat((0,0,0),(1,0,0),(0,0,0));
14 l103:=we*mat((0,0,0),(0,0,0),(1,0,0));
15 l112:=we*mat((0,0,0),(1,0,0),(0,0,0));
16 l121:=we*mat((0,0,0),(0,0,0),(1,0,0));
17
18 % nonlinear M matrices
19 mm000:=mat((0,0,0),(0,0,0),(0,0,0));
20 mm001:=mm000;
21 mm002:=mm000;
22 mm100:=mat((0,-1,0),(-1,0,0),(0,0,0));
23 mm101:=mat((0,0,0),(0,-1,0),(0,0,0));
24 mm102:=mat((0,0,0),(0,0,0),(0,-1,0));
25 mm010:=mat((0,0,-1),(0,0,0),(-1,0,0));
26 mm011:=mat((0,0,0),(0,0,-1),(0,0,0));
27 mm012:=mat((0,0,0),(0,0,0),(0,0,-1));
```

Now specify the eigenvalues of $\mathfrak{L}_{\vec{0}}$, thus defining $A_{\vec{0}}$ (Definition 4), the single eigenvalue (since $\mathfrak{m}=1$) of the slow subspace, and $B_{\vec{0}}$, the $(d-1)\times(d-1)$ matrix of eigenvalues of the stable subspace (from Definition 4), $\mathfrak{m}'=d-1$).

```
28 array aa(nn,nn);
29 aa(0,0):=0$
30 dr:=d-1$
31 array bb(dr-1,dr-1,nn,nn);
32 bb(0,0,0,0):=-b+i*f$
33 bb(1,1,0,0):=-b-i*f$
```

For coding purposes, it is more convenient to define the matrices $\mathfrak{L}_{\vec{n}}$ and $\mathfrak{M}_{\vec{n},j}$ as arrays.

```
34 kmax:=3$ % highest order spatial derivative
35 array lls(dr,dr,kmax,kmax);
36 for k1:=0:dr do for k2:=0:dr do lls(k1,k2,0,0):=ll00(k1+1,k2+1);
37 for k1:=0:dr do for k2:=0:dr do lls(k1,k2,1,0):=ll10(k1+1,k2+1);
38 for k1:=0:dr do for k2:=0:dr do lls(k1,k2,0,1):=ll01(k1+1,k2+1);
39 for k1:=0:dr do for k2:=0:dr do lls(k1,k2,2,0):=ll20(k1+1,k2+1);
40 for k1:=0:dr do for k2:=0:dr do lls(k1,k2,0,2):=ll02(k1+1,k2+1);
41 for k1:=0:dr do for k2:=0:dr do lls(k1,k2,3,0):=ll30(k1+1,k2+1);
42 for k1:=0:dr do for k2:=0:dr do lls(k1,k2,0,3):=ll03(k1+1,k2+1);
43 for k1:=0:dr do for k2:=0:dr do lls(k1,k2,1,2):=ll12(k1+1,k2+1);
44 for k1:=0:dr do for k2:=0:dr do lls(k1,k2,2,1):=ll21(k1+1,k2+1);
45 array mms(dr,dr,dr,kmax,kmax);
46 for k1:=0:dr do for k2:=0:dr do mms(0,k1,k2,0,0):=mm000(k1+1,k2+1);
47 for k1:=0:dr do for k2:=0:dr do mms(0,k1,k2,1,0):=mm100(k1+1,k2+1);
48 for k1:=0:dr do for k2:=0:dr do mms(0,k1,k2,0,1):=mm010(k1+1,k2+1);
49 for k1:=0:dr do for k2:=0:dr do mms(1,k1,k2,0,0):=mm001(k1+1,k2+1);
50 for k1:=0:dr do for k2:=0:dr do mms(1,k1,k2,1,0):=mm101(k1+1,k2+1);
51 for k1:=0:dr do for k2:=0:dr do mms(1,k1,k2,0,1):=mm011(k1+1,k2+1);
52 for k1:=0:dr do for k2:=0:dr do mms(2,k1,k2,0,0):=mm002(k1+1,k2+1);
53 for k1:=0:dr do for k2:=0:dr do mms(2,k1,k2,1,0):=mm102(k1+1,k2+1);
54 for k1:=0:dr do for k2:=0:dr do mms(2,k1,k2,0,1):=mm012(k1+1,k2+1);
```

## A.2 Generalised eigenvectors and slow PDE

The code below derives the centre (or slow) and stable modes for some PDE of the form (6), and constructs the slow PDE projected onto the centre (or slow) subspace (as describe in Section 3.2). The code is applicable to any $d$ dimensional microscale field on an $M = 2$ dimensional 'large' domain with a PDE of nonlinear order $P = 2$, where the $d$ dimensional field separates into a $m' = d - 1$ stable subspace and a $m = 1$ dimensional slow or centre subspace (although v, defined in line 63, is list of at least $d$ dummy variables, so the number of dummy variables in this list should be increased if field u has dimension $d > 7$, and similarly for w defined in line 97, is list of at least $d(d-1)$ dummy variables).

This code is used to derive the slow PDE of the fluid flow problem described in Section 3.2, with parameters defined in Appendix A.1

Calculate the centre (or slow) left and right eigenvectors, $Z^{\vec{0}}$ and $V^{\vec{0}}$, of $\mathfrak{L}_{\vec{0}}$ from the nullspaces of $[\mathfrak{L}_{\vec{0}} - A_{\vec{0}}]^{\mathsf{T}}$ and $[\mathfrak{L}_{\vec{0}} - A_{\vec{0}}]$, respectively (Definition 4).

```
55 matrix ii(d,d)$ % identity matrix
56 for j:=1:d do ii(j,j):=1;
57 vvec:=first(nullspace (ll00-ii*aa(0,0)))$ % right eigenvec
58 zvec:=first(nullspace tp(ll00-ii*aa(0,0)))$ % left eigenvec
59 zvec:=zvec/det(tp(zvec)*vvec)$ % rescale
60 array vvt(dr,nn,nn);
61 for j:=0:dr do vvt(j,0,0):=vvec(j+1,1);
62 array zz0(dr);
63 for j:=0:dr do zz0(j):=zvec(j+1,1);
64 v:={v0, v1, v2, v3, v4, v5, v6, v7}$ % dummy variables
```

The recurrence relation (28) is applied to construct the slow eigenvectors $\tilde{\mathcal{V}}^{\vec{n}}$ and eigenvalues $A_{\vec{n}}$ of $\tilde{\mathcal{L}}$ for all $0 < |\vec{n}| \leqslant N$. Once we know all $A_{\vec{n}}$, we know the linear part of the slow PDE (44) (or linear parts of the slow PDEs (4a) and (5a) for the thin film flow example).

```
65 for p1:=0:nn do for p2:=0:(nn-p1) do if (p1+p2)>0 then <<
66     aa(p1,p2):=(for k1:=0:min(kmax,p1) sum for k2:=0:min(kmax,p2) sum
67         for j1:=0:dr sum for j2:=0:dr sum
68         zz0(j1)*lls(j1,j2,k1,k2)*vvt(j2,p1-k1,p2-k2));
69     aa(p1,p2):=(aa(p1,p2) where {xi1=>0, xi2=>0});
70     eqn:={(for j1:=0:dr sum
71         zz0(j1)*part(v,j1+1))-xi1^p1*xi2^p2
72         /factorial(p1)/factorial(p2)};
73     for j1:=0:dr do <<
74         co:=(for j2:=0:dr sum lls(j1,j2,0,0)*part(v,j2+1) )
75             -part(v,j1+1)*aa(0,0)
76             +(for k1:=0:min(kmax,p1) sum
77             for k2:=0:min(kmax,p2) sum for j2:=0:dr sum
78             lls(j1,j2,k1,k2)*vvt(j2,p1-k1,p2-k2))
79             -(for k1:=0:p1 sum for k2:=0:p2 sum
80             vvt(j1,p1-k1,p2-k2)*aa(k1,k2));
81         eqn:=co.eqn;
82     >>;
83     solv:=solve(eqn,v);
84     for j1:=0:dr do vvt(j1,p1,p2):=sub(solv,part(v,j1+1));
85 >>;
```

Now consider the stable subspace. As was done above for slow subspace, first calculate the left and right eigenvectors associated with the eigenvalues in $B_{\vec{0}}$ of matrix $\mathfrak{L}_{\vec{0}}$, and then apply a recurrence relation analogous to (28) to construct the stable eigenvectors $\tilde{\mathcal{W}}^{\vec{n}}$ and associated matrices of eigenvalues $B_{\vec{n}}$ of $\tilde{\mathcal{L}}$ for all $0 < |\vec{n}| \leqslant N$.

```
 86 array wwt(dr,dr-1,nn,nn);
 87 for j1:=0:(dr-1) do <<
 88     wvec:=first(nullspace (ll00-ii*bb(j1,j1,0,0))); % right eigenvec
 89     for j2:=0:dr do wwt(j2,j1,0,0):=wvec(j2+1,1);
 90 >>;
 91 array ww0t(dr-1,dr);
 92 for j1:=0:(dr-1) do <<
 93     wvec:=first(nullspace tp(ll00-ii*bb(j1,j1,0,0))); % left eigenvec
 94     scl:=(for j2:=0:dr sum wvec(j2+1,1)*wwt(j2,j1,0,0));
 95     for j2:=0:dr do ww0t(j1,j2):=wvec(j2+1,1)/scl;
 96 >>;
 97 w:={w0, w1, w2, w3, w4, w5, w6, w7, w8, w9, w10, w11, w12, w13, w14}$
 98 array wtmp(dr,dr-1);
 99 for j1:=0:dr do for j2:=0:(dr-1) do wtmp(j1,j2):=part(w,j1*dr+j2+1);
100
101 operator dirac;
102 let {dirac(~j1,~j2)=>0 when j1 neq j2, dirac(~j1,~j2)=>1 when j1=j2};
103 for p1:=0:nn do for p2:=0:(nn-p1) do if (p1+p2)>0 then <<
104     for m1:=0:(dr-1) do for m2:=0:(dr-1) do bb(m1,m2,p1,p2):=
105     (for k1:=0:min(kmax,p1) sum for k2:=0:min(kmax,p2) sum
106     for j1:=0:dr sum for j2:=0:dr sum
107     ww0t(m1,j1)*lls(j1,j2,k1,k2)*wwt(j2,m2,p1-k1,p2-k2));
108     for m1:=0:(dr-1) do for m2:=0:(dr-1) do
109     bb(m1,m2,p1,p2):=(bb(m1,m2,p1,p2) where {xi1=>0, xi2=>0});
110     eqn:={};
111     for m1:=0:(dr-1) do for m2:=0:(dr-1) do <<
112         co:=(for j1:=0:dr sum
113             ww0t(m1,j1)*wtmp(j1,m2))-
114             dirac(m1,m2)*xi1^p1*xi2^p2/factorial(p1)/factorial(p2);
115         eqn:=co.eqn;
116     >>;
117     for j1:=0:dr do for j3:=0:(dr-1) do <<
118         co:=(for j2:=0:dr sum lls(j1,j2,0,0)*wtmp(j2,j3))
119             -(for j2:=0:(dr-1) sum wtmp(j1,j2)*bb(j2,j3,0,0))
120             +(for k1:=0:min(kmax,p1) sum
121             for k2:=0:min(kmax,p2) sum for j2:=0:dr sum
122             lls(j1,j2,k1,k2)*wwt(j2,j3,p1-k1,p2-k2))
123             -(for k1:=0:p1 sum for k2:=0:p2 sum for j2:=0:(dr-1) sum
124             wwt(j1,j2,p1-k1,p2-k2)*bb(j2,j3,k1,k2));
125         eqn:=co.eqn;
126     >>;
127     solw:=solve(eqn,w);
128     for j1:=0:dr do for j2:=0:(dr-1) do wwt(j1,j2,p1,p2)
129         :=sub(solw,part(w,j1*dr+j2+1));
130 >>;
```

Construct the nonlinear term $\tilde{f}[\tilde{u}]$ (25) with $\tilde{u} = \mathcal{V}\mathcal{U} + \mathcal{W}\mathcal{S}$ as in (30), where $\mathcal{V} = [\mathcal{V}^{\vec{n}}]$, $\mathcal{W} = [\mathcal{W}^{\vec{n}}]$, $\mathcal{U} = [U_{\vec{n}}]$, and $\mathcal{S} = [S_{\vec{n}}]$. Here we ignore the uncertain gradient terms (second line of (25)) because such terms are regulated to the remainder $\rho$ (45) of the slow PDE and because we have chosen N large enough so that these uncertain terms only appear in the remainder's convolution term. Then, compare coefficients of multinomial variable $\vec{\xi}$ to determine the slow and stable parts of the nonlinear term, $f_c$ and $f_s$, as in (31).

```
131 operator uu;
132 operator ss;
133 operator nonl;
134 operator nonlvw;
135 operator tu;
136 for j:=0:dr do tu(j):=(for k1:=0:nn sum for k2:=0:(nn-k1) sum
137     (vvt(j,k1,k2)*uu(k1,k2)+(for n1:=0:(dr-1) sum
138         wwt(j,n1,k1,k2)*ss(n1,k1,k2))) where {b^2=>bf-f^2});
139 for l:=0:dr do nonl(l):=(for p1:=0:min(kmax,nn) sum
140     for p2:=0:min(kmax,nn) sum
141     for n1:=0:nn sum for n2:=0:(nn-n1) sum
142     for m1:=0:n1 sum for m2:=0:n2 sum
143         xi1^(n1)*xi2^(n2)/(factorial(m1)*factorial(m2)
144         *factorial(n1-m1)*factorial(n2-m2))
145         *(for i:=0:dr sum for j:=0:dr sum
146         sub(xi1=0,xi2=0,df(df(tu(i),xi1,m1+p1),xi2,m2+p2))
147         *mms(l,i,j,p1,p2)*sub(xi1=0,xi2=0,tu(j)) ) );
148 for l:=0:dr do nonlvw(l):=(for k1:=0:nn sum for k2:=0:(nn-k1) sum
149     vvt(l,k1,k2)*fc(k1,k2)+
150     (for j:=0:(dr-1) sum wwt(l,j,k1,k2)*fs(j,k1,k2)  ));
151
152 operator fc;
153 operator fs;
154 for l:=0:dr do nonlvw(l):=(for k1:=0:nn sum for k2:=0:(nn-k1) sum
155     vvt(l,k1,k2)*fc(k1,k2)+
156     (for j:=0:(dr-1) sum wwt(l,j,k1,k2)*fs(j,k1,k2)  ));
157 ffs:={}$
158 for k1:=0:nn do for k2:=0:(nn-k1) do ffs:=fc(k1,k2).ffs;
159 for k1:=0:nn do for k2:=0:(nn-k1) do for j:=0:(dr-1) do
160     ffs:=fs(j,k1,k2).ffs;
161 eqn:={}$
162 for l:=0:dr do for k1:=0:nn do for k2:=0:(nn-k1) do
163     eqn:=coeffn(coeffn(nonl(l)-nonlvw(l),xi1,k1),xi2,k2).eqn;
164 eqn:=(eqn where b^2=>bf-f^2)$
165 solf:=solve(eqn,ffs)$

166 operator uu;
167 operator ss;
168 operator nonl;
169 operator nonlvw;
170 operator tu;
171 for j:=0:dr do tu(j):=(for k1:=0:nn sum for k2:=0:(nn-k1) sum
172     (vvt(j,k1,k2)*uu(k1,k2)+(for n1:=0:(dr-1) sum
173         wwt(j,n1,k1,k2)*ss(n1,k1,k2))) where {b^2=>bf-f^2});
174 for l:=0:dr do nonl(l):=(for p1:=0:min(kmax,nn) sum
175     for p2:=0:min(kmax,nn) sum
176     for n1:=0:nn sum for n2:=0:(nn-n1) sum
```

```
177     for m1:=0:n1 sum for m2:=0:n2 sum
178         xi1^(n1)*xi2^(n2)/(factorial(m1)*factorial(m2)
179         *factorial(n1-m1)*factorial(n2-m2))
180         *(for i:=0:dr sum for j:=0:dr sum
181         sub(xi1=0,xi2=0,df(df(tu(i),xi1,m1+p1),xi2,m2+p2))
182         *mms(l,i,j,p1,p2)*sub(xi1=0,xi2=0,tu(j)) ) );
183
184 operator fc;
185 operator fs;
186 operator ftmp;
187 ffs:={ftmp(0)}$
188 for j:=1:dr do ffs:=ftmp(j).ffs;
189 for k1:=nn step -1 until 0 do for k2:=nn-k1 step -1 until 0 do <<
190 eqn:={}$
191 for l:=0:dr do nonlvw(l):=coeffn(coeffn((vvt(l,k1,k2)*ftmp(0)+
192             (for j:=0:(dr-1) sum wwt(l,j,k1,k2)
193                 *ftmp(j+1))),xi1,k1),xi2,k2);
194 for l:=0:dr do eqn:=coeffn(coeffn(nonl(l)-nonlvw(l)
195     ,xi1,k1),xi2,k2).eqn;
196 solf:=solve(eqn,ffs);
197 fc(k1,k2):=sub(ffs,ftmp(0));
198 for l:=0:(dr-1) do fs(l,k1,k2):=sub(ffs,ftmp(l+1));
199 >>;
200
201
```

For the general slow PDE (44) of $U(\vec{x}, t)$ (or example PDEs (4a) and (5a)) we require the projection onto the centre (or slow) subspace $\mathbb{E}_c^0$. In addition, as discussed in Section 3.2, we regulate all dependence on the stable subspace to a remainder $\rho$ (45), which, for deriving the slow PDE effectively requires us to set parameter $S$ to zero. In the code below we project nonlinear terms $f_c(\vec{x}, t)$ and $f_s(\vec{x}, t)$ onto $\mathbb{E}_c^0$ and set $S$ to zero, as defined in equation (40). Then, using $U_{\vec{n}} = \partial_{\vec{x}}^{\vec{n}} U(\vec{x}, t)$, from (38) but neglecting the convolutions (which are also placed in remainder $\rho$ (45)), we rewrite the remaining nonlinear terms $f_c^{\vec{0}}(\vec{x}, t) + f_s(\vec{x}, t)$ in the form $a_{\vec{n},\vec{m}} \partial_{\vec{x}}^{\vec{n}} U(\vec{x}, t) \partial_{\vec{x}}^{\vec{m}} U(\vec{x}, t)$ (43).

```
202 fc0u:=(sub(solf,fc(0,0)) where {ss(~j1,~j2,~j3)=>0, bf=>b^2+f^2})$
203 fs0u:=(for k1:=0:nn sum for k2:=0:(nn-k1) sum
204     for j:=0:dr sum for l:=0:(dr-1) sum zz0(j)
205     *wwt(j,l,k1,k2)*sub(solf,fs(l,k1,k2))
206     where {xi1=>0, xi2=>0, ss(~j1,~j2,~j3)=>0, bf=>b^2+f^2})$
207
208 fctmp:=fc0u+fs0u$
209 array ac(nn,nn,nn,nn);
210 for k1:=0:nn do for k2:=0:(nn-k1) do
211     for j1:=0:(k1+k2) do for j2:=0:(k1+k2-j1) do <<
212         if (k1+k2=j1+j2 and ac(j1,j2,k1,k2)=0) then
213             ac(k1,k2,j1,j2):=coeffn(coeffn(fctmp,uu(j1,j2),1)
214                 ,uu(k1,k2),1)/2
215         else ac(k1,k2,j1,j2):=coeffn(coeffn(fctmp,uu(j1,j2),1)
216             ,uu(k1,k2),1);
217         fctmp:=fctmp-ac(k1,k2,j1,j2)*uu(j1,j2)*uu(k1,k2);
218 >>;
```

```
219 for k1:=0:nn do for k2:=0:(nn-k1) do <<
220     ac(k1,k2,k1,k2):=coeffn(fctmp,uu(k1,k2),2);
221     fctmp:=fctmp-ac(k1,k2,k1,k2)*uu(k1,k2)^2;
222 >>;
```

We now write all $A_{\vec{n}}$ and all elements of $\mathcal{C}$ to define the slow PDE (44) of $U(\vec{x},t)$ (or example PDEs (4a) and (5a)).

```
223 for k1:=0:nn do for k2:=0:(nn-k1) do
224     if aa(k1,k2) neq 0 then write aa(k1,k2):=aa(k1,k2);
225 for k1:=0:nn do for k2:=0:(nn-k1) do
226     for j1:=0:nn do for j2:=0:(nn-j1) do
227     if ac(k1,k2,j1,j2) neq 0 then write ac(k1,k2,j1,j2)
228         :=ac(k1,k2,j1,j2);
```

## A.3  Remainder of the slow PDE

In this section we calculate some terms which appear in the remainder $\rho$ (45) to illustrate how this remainder is constructed for a specific model. We assume that the remainder is evaluated after some relatively long time so that terms $\mathcal{O}(e^{-\gamma t})$ are negligible; as $\gamma \in (\alpha, \beta)$, this ensures that the linear effects of the stable modes are negligible (Assumption 3).

We calculate all inner products $\langle Z^{\vec{0}}, \tilde{\mathcal{W}}^{\vec{n}} \rangle_{\vec{\xi}=\vec{0}}$ and $\langle Z^{\vec{0}}, \partial_{\vec{\xi}}^{\vec{n}} \tilde{\mathcal{W}}^{\vec{k}} \rangle_{\vec{\xi}=\vec{0}}$ and use these to calculate

$$\sum_{|\vec{n}|=0}^{N} A_{\vec{n}} \langle Z^{\vec{0}}, \partial_{\vec{\xi}}^{\vec{n}} \tilde{\mathcal{W}}^{\vec{k}} \rangle_{\vec{\xi}=\vec{0}} \quad \text{and} \quad \left[ \langle Z^{\vec{0}}, \tilde{\mathcal{W}} \rangle_{\vec{\xi}=\vec{0}} \mathcal{B} \right]_{\vec{k}} = \sum_{|\vec{n}|=0}^{N} \langle Z^{\vec{0}}, \tilde{\mathcal{W}}^{\vec{n}} \rangle_{\vec{\xi}=\vec{0}} B_{\vec{k}-\vec{n}} \,.$$

```
229 array ipzw(dr-1,nn,nn);
230 for k1:=0:nn do for k2:=0:(nn-k1) do
231     for l:=0:(dr-1) do ipzw(l,k1,k2):=(for j:=0:dr sum
232         zz0(j)*wwt(j,l,k1,k2) where {xi1=>0, xi2=>0})$
233 
234 array ipzdw(dr-1,nn,nn,nn,nn);
235 array asum(dr-1,nn,nn);
236 for k1:=0:nn do for k2:=0:(nn-k1) do for j1:=0:k1 do for j2:=0:k2 do
237     for l:=0:(dr-1) do ipzdw(l,k1,k2,j1,j2):=
238     (for j:=0:dr sum zz0(j)*df(df(wwt(j,l,k1,k2),xi1,j1),xi2,j2))$
239 for k1:=0:nn do for k2:=0:(nn-k1) do
240     for l:=0:(dr-1) do asum(l,k1,k2):=(for j1:=0:nn sum
241     for j2:=0:(nn-j1) sum aa(j1,j2)*ipzdw(l,k1,k2,j1,j2)
242     where {xi1=>0, xi2=>0})$
243 
244 array ipzbb(dr-1,nn,nn);
245 for k1:=0:nn do for k2:=0:(nn-k1) do
246     for l:=0:(dr-1) do <<
247         ipzbb(l,k1,k2):=(for j:=0:(dr-1) sum for j1:=0:k1 sum
248         for j2:=0:k2 sum
249         ipzw(j,j1,j2)*bb(j,l,k1-j1,k2-j2)
250             where {xi1=>0, xi2=>0})$
251     if  ipzbb(l,k1,k2) neq 0 then write {k1,k2};
252 >>;
```

For example, $\langle Z^{\vec{0}}, \tilde{\mathcal{W}}^{(1,2)} \rangle_{\vec{\xi}=\vec{0}}$ is

```
253 mat((ipzw(0,1,2)),(ipzw(1,1,2)));
```
$\sum_{|\vec{n}|=0}^{N} A_{\vec{n}} \langle Z^{\vec{0}}, \partial_{\vec{\xi}}^{\vec{n}} \tilde{\mathcal{W}}^{(1,2)} \rangle_{\vec{\xi}=\vec{0}}$ is
```
254 mat((asum(0,1,2)),(asum(1,1,2)));
```
and $[\langle Z^{\vec{0}}, \tilde{\mathcal{W}} \rangle_{\vec{\xi}=\vec{0}} \mathcal{B}]_{(10)}$ and $[\langle Z^{\vec{0}}, \tilde{\mathcal{W}} \rangle_{\vec{\xi}=\vec{0}} \mathcal{B}]_{(01)}$ are
```
255 mat((ipzbb(0,1,0)),(ipzbb(1,1,0)));
256 mat((ipzbb(0,0,1)),(ipzbb(1,0,1)));
```
The last line of (45), in terms of $\mathcal{S}$ (which contains the convolution (34)) is
```
257 convterm:=for k1:=0:nn sum for k2:=0:(nn-k1) sum
258     for l:=0:(dr-1) sum
259     ((for n1:=0:k1 sum for n2:=0:k2 sum for j:=0:(dr-1) sum
260         ipzw(j,n1,n2)*bb(j,l,k1-n1,k2-n2) )
261         -asum(l,k1,k2))*ss(l,k1,k2);
```

Finally we calculate linear uncertain term $r_c^{\vec{0}}(t) + \langle Z^{\vec{0}}, \tilde{\mathcal{W}} \rangle_{\vec{\xi}=\vec{0}} r_s(t)$ and show that it is zero for our chosen order $N > \max(p_i^j, k_{\max}) = 1$. To do this we first construct the forcing $\tilde{r}[u]$ (22) in terms of the uncertain gradient terms, that is, spatial derivatives of $u^{\vec{n}}$ with $|\vec{n}| = N$, and compare this to $\tilde{r}[u] = \tilde{\mathcal{V}} r_c(t) + \tilde{\mathcal{W}} r_s(t)$ (30) to determine all components of $r_c(t)$ and $r_s(t)$.

```
262 operator ux;
263 depend ux, x1, x2;
264 array ru(dr);
265 for k1:=0:kmax do for k2:=0:(kmax-k1) do for j1:=0:dr do
266     for j2:=0:dr do
267     for n1:=0:nn do for n2:=0:(nn-n1) do for l1:=0:nn do
268     if (l1 leq n1+k1 and nn-l1 leq n2+k2 and
269         n1+n2+k1+k2 > nn and k1+k2 > 0) then
270         ru(j1):=ru(j1)+lls(j1,j2,k1,k2)*xi1^n1*xi2^n2/
271         (factorial(n1)*factorial(n2))*factorial(k1+n1)*factorial(k2+n2)
272         /(factorial(l1)*factorial(nn-l1))
273         *df(df(ux(l1,nn-l1,j2),x1,k1+n1-l1),x2,k2+n2-nn+l1);
274
275 operator rvw;
276 operator rc;
277 operator rs;
278 for l:=0:dr do rvw(l):=(for k1:=0:nn sum
279     for k2:=0:(nn-k1) sum
280     vvt(l,k1,k2)*rc(k1,k2)+
281     (for j:=0:(dr-1) sum wwt(l,j,k1,k2)*rs(j,k1,k2) ));
282 rrs:={}$
283 for k1:=0:nn do for k2:=0:(nn-k1) do rrs:=rc(k1,k2).rrs;
284 for k1:=0:nn do for k2:=0:(nn-k1) do for j:=0:(dr-1) do
285     rrs:=rs(j,k1,k2).rrs;
286 eqn:={}$
287 for l:=0:dr do for k1:=0:nn do for k2:=0:(nn-k1) do
288     eqn:=coeffn(coeffn(ru(l)-rvw(l),xi1,k1),xi2,k2).eqn;
289 eqn:=(eqn where b^2=>bf-f^2)$
290 solr:=solve(eqn,rrs)$
291 solr:=(solr where bf=>b^2+f^2)$
```

In the remainder $\rho$ we only retain the $r_c^{\vec{0}}$ component of the centre subspace forcing $r_c(t)$, and for the stable subspace forcing we calculate $\langle Z^{\vec{0}}, \tilde{\mathcal{W}} \rangle_{\vec{\xi}=\vec{0}} r_s(t) = \sum_{|\vec{n}|=0}^{N} \langle Z^{\vec{0}}, \tilde{\mathcal{W}}^{\vec{n}} \rangle r_s^{\vec{n}}(t)$. Then we show that the sum $r_c^{\vec{0}} + \langle Z^{\vec{0}}, \tilde{\mathcal{W}} \rangle_{\vec{\xi}=\vec{0}} r_s(t) = 0$.

```
292 rc0:=sub(solr,rc(0,0));
```

```
293 rsall:=(for j:=0:(dr-1) sum for k1:=0:nn sum
294     for k2:=0:(nn-k1) sum
295     ipzw(j,k1,k2)*sub(solr,rs(j,k1,k2)));
296 zero:=rc0+rsall;
297
298 end;
```

Here we do not calculate the nonlinear terms which appear in the remainder $\rho$ (45), but these terms can be calculated by editing the nonlinear calculations in Appendix A.2 to include the second line of (24). As before, equation (31) separates the centre and stable components of $\tilde{f}[\tilde{u}]$, and (40) separates the nonlinear components into terms which are included in the slow PDE (44) and those which are in the remainder $\rho$ (45). As the nonlinear calculations in Appendix A.2 are already fairly memory intensive, we do not extend these calculations to deriving the nonlinear remainder terms in $\rho$.